\documentclass[11pt]{amsart}

\usepackage[latin1]{inputenc}
\usepackage{amsmath, latexsym, amssymb}
\usepackage{graphicx}
\usepackage{subfigure}
\usepackage{natbib}
\usepackage{color}
\usepackage{amsmath, amsfonts, amssymb, hyperref}
\usepackage{graphics}
\newcommand{\Q}{\mathcal{Q}} 
\definecolor{gr75}{gray}{0.75}
\newcommand{\gx}[1]{\mbox{\textbf{#1}}}
 
\numberwithin{equation}{section}
\newcommand{\augdg}{\widehat{dg}}

\newcommand{\ms}{\begin{math}}
\newcommand{\me}{\end{math}}
\newcommand{\da}{\mathcal{A}}
\newcommand{\qs}{\mathcal{S}}
\newcommand{\fo}{\mathcal{F}o}
\newcommand{\btr}{\blacktriangleright}
\newcommand{\D}{\mathcal{D}}

\newlength\cellsize 
\savebox2{%
\begin{picture}(15,15)
\put(0,0){\line(1,0){15}}
\put(0,0){\line(0,1){15}}
\put(15,0){\line(0,1){15}}
\put(0,15){\line(1,0){15}}
\end{picture}}
\newcommand\cellify[1]{\def\thearg{#1}\def\nothing{}%
\ifx\thearg\nothing
\vrule width0pt height\cellsize depth0pt\else
\hbox to 0pt{\usebox2\hss}\fi%
\vbox to 15\unitlength{
\vss
\hbox to 15\unitlength{\hss$#1$\hss}
\vss}}
\newcommand\tableau[1]{\vtop{\let\\=\cr
\setlength\baselineskip{-16000pt}
\setlength\lineskiplimit{16000pt}
\setlength\lineskip{0pt}
\halign{&\cellify{##}\cr#1\crcr}}}
\savebox3{%
\begin{picture}(15,15)
\put(0,0){\line(1,0){15}}
\put(0,0){\line(0,1){15}}
\put(15,0){\line(0,1){15}}
\put(0,15){\line(1,0){15}}
\end{picture}}
\newcommand\expath[1]{%
\hbox to 0pt{\usebox3\hss}%
\vbox to 15\unitlength{
\vss
\hbox to 15\unitlength{\hss$#1$\hss}
\vss}}

\theoremstyle{plain} 
\newtheorem{theorem}{Theorem}[section]
\newtheorem{proposition}[theorem]{Proposition}

\newtheorem{lemma}[theorem]{Lemma}
\newtheorem{corollary}[theorem]{Corollary}

\newtheorem{definition}[theorem]{Definition}
\newtheorem*{example}{Example}

\newtheorem*{remark}{Remark}

\begin{document}
\title[Quasisymmetric Schur functions]{Quasisymmetric Schur functions}

\author{J. Haglund}
\address{Department of Mathematics, University of Pennsylvania, Philadelphia, PA 19104-6395, USA}
\email{\href{jhaglund@math.upenn.edu}{jhaglund@math.upenn.edu}}

\author{K. Luoto}
\address{Department of Mathematics, University of British Columbia, Vancouver, BC V6T 1Z2, Canada}
\email{\href{mailto:kwluoto@math.ubc.ca}{kwluoto@math.ubc.ca}}

\author{S. Mason}
\address{Department of Mathematics, Davidson College, Davidson, NC 28035-7129, USA}
\email{\href{sarahkmason@gmail.com }{sarahkmason@gmail.com}}

\author{S. van Willigenburg}
\address{Department of Mathematics, University of British Columbia, Vancouver, BC V6T 1Z2, Canada}
\email{\href{mailto:steph@math.ubc.ca}{steph@math.ubc.ca}}
\thanks{The first author was supported in part by  NSF grants DMS 0553619 and DMS 0901467. The third author was supported in part by NSF postdoctoral research fellowship DMS 0603351. The second and fourth authors were supported in part by the National Sciences and Engineering Research Council of Canada. The authors would like to thank the Banff International Research Station and the Centre de Recherches Math\'{e}matiques, where some of the research took place.}
\subjclass[2000]{Primary 05E05; Secondary 05E10, 33D52} 
\keywords{compositions, Kostka coefficients, nonsymmetric Macdonald polynomials, Pieri rule, quasisymmetric function, Schur function, tableaux} 
\begin{abstract}
We introduce a new basis for quasisymmetric functions, which arise from a specialization of nonsymmetric Macdonald polynomials to standard bases, also known as Demazure atoms. Our new basis is called the basis of quasisymmetric Schur functions, since the basis elements refine Schur functions in a natural way. We derive expansions for quasisymmetric Schur functions in terms of monomial and fundamental quasisymmetric functions, which give rise to quasisymmetric refinements of Kostka numbers and standard (reverse) tableaux. From here we derive a Pieri rule for quasisymmetric Schur functions that naturally refines the Pieri rule for Schur functions.  After surveying 
combinatorial formulas for Macdonald polynomials, including an
expansion
of Macdonald polynomials into fundamental quasisymmetric functions, we show how some
of our results can be extended to include the
$t$ parameter from Hall-Littlewood theory.
\end{abstract}

\maketitle
\tableofcontents
\section{Introduction}\label{sec:intro} 
Macdonald polynomials were originally introduced in 1988 \cite{Mac88, Mac}, as a solution to a problem involving Selberg's integral posed by  Kadell \cite{Kad88}. They are $q,t$ analogues of symmetric functions such that setting $q=t=0$ in the Macdonald polynomial $P_\lambda (X; q, t)$, for $\lambda$ a partition, yields the Schur function $s_\lambda$. Since their introduction they have arisen in further mathematical areas such as representation theory and quantum computation. For example,  Cherednik \cite{Cherednik} showed that  nonsymmetric Macdonald polynomials are connected to the representation theory of double affine Hecke algebras, and setting $q=t^\alpha$, dividing by a power of $1-t$ and letting $t \to 1$ yields Jack polynomials, which model bosonic variants of single component abelian and nonabelian fractional quantum Hall states \cite{BH}. The aforementioned nonsymmetric Macdonald polynomials, $E ^{\prime}_\alpha (X; q,t)$ where $\alpha$ is a weak composition, are  a nonsymmetric refinement of the $P_\lambda (X; q, t)$. Setting $q=t=0$ in an identity of Macdonald and Marshall expressing $P_\lambda (X; q, t)$ as a linear combination
of modified versions of the $E^{\prime}$'s (see Section \ref{sec:conc})
implies that Schur functions can be decomposed into nonsymmetric functions $\da _\gamma$ for $\gamma$ a weak composition. These functions were first studied in \cite{LS}, where they were termed standard bases, however, to avoid confusion with other objects termed standard bases, we refer to them here as Demazure atoms since they decompose Demazure characters into their smallest parts. The definition we use also differs from that in \cite{LS} as our definition not only is arguably simpler than the one appearing there, but also is upward compatible with the new combinatorics appearing in the combinatorial formulae for Type $A$  symmetric and nonsymmetric Macdonald polynomials \cite{HHL2, HHL}. The equivalence of these two definitions is established in \cite{Mason2}. It should be stressed that Demazure \emph{atoms} should not be confused with Demazure \emph{characters}, which involve the combinatorial tool of crystal graphs.  However, certain linear combinations of Demazure atoms  form Demazure characters, and their relationship to each other and to nonsymmetric Macdonald polynomials can be found in \cite{Ion, Ion2, Mason2}.

Interpolating between symmetric functions and nonsymmetric functions are quasisymmetric functions. These were introduced  as a source of generating functions for $P$-partitions \cite{Gessel} but since then, like Macdonald polynomials, they have impacted, and deepened the understanding of, other areas. For example in category theory they are a terminal object in the category of graded Hopf algebras equipped with a zeta function \cite{Aguiar}; in lattice theory they induce Pieri rules analogous to those found in the algebra of symmetric functions \cite{BMSvW}; in discrete geometry the quasisymmetric functions known as peak functions were found to be dual to the \textbf{cd}-index \cite{BHvW}; in symmetric function theory they identify equal ribbon Schur functions \cite{HDL}; in representation theory they arise as characters of a degenerate quantum group \cite{Hivert, KrobThibon}.

Therefore, a natural object to seek is a quasisymmetric function that interpolates between the nonsymmetric Schur functions, known as Demazure atoms, and Schur functions. Furthermore, since Demazure atoms exhibit many Schur function properties \cite{Mason1}, a natural question to ask is which properties of Schur functions are exhibited by \emph{quasisymmetric} Schur functions? In this paper we define quasisymmetric Schur functions and show they naturally lift well known combinatorial properties of symmetric functions indexed by partitions, to combinatorial properties of quasisymmetric functions indexed by compositions. More precisely, we show the following.
\begin{enumerate}
\item The expression for Schur functions in terms of monomial symmetric functions refines to an expression for quasisymmetric Schur functions in terms of monomial quasisymmetric functions, giving rise to quasisymmetric Kostka coefficients.
\item The expression for Schur functions in terms of fundamental quasisymmetric functions naturally refines to quasisymmetric Schur functions.
\item The Pieri rule for multiplying a Schur function indexed by a row or a column with a generic Schur function refines to a rule for multiplying a quasisymmetric Schur function indexed by a row or a column with a generic quasisymmetric Schur function. Moreover, this rule is a new example of the construction studied in \cite{BMSvW, BHvW}, where the underlying poset involved is a poset of compositions.
\end{enumerate}
The existence of such results introduces a plethora of research avenues to pursue concerning the quasisymmetric analogues of other symmetric function properties.
For example, the latter result naturally raises the question of whether the Littlewood-Richardson rule for multiplying two generic Schur functions can be refined to quasisymmetric Schur functions. Such a refinement may not be easy to find as the classical Littlewood-Richardson rule produces nonnegative structure constants, whereas multiplying together two quasisymmetric Schur functions sometimes results in negative structure constants. The smallest such example exists at $n=6$. However, in the sequel to this paper we successfully refine the Littlewood-Richardson rule by multiplying a generic Schur function and quasisymmetric Schur function \cite{HLMvW2}.

More precisely, this paper is structured as follows. In Sections~\ref{sec:qsymsym}, \ref{sec:schur}, \ref{sec:atoms} we review the necessary, and sometimes nonstandard, background material regarding quasisymmetric and symmetric functions, and Demazure atoms. In Section~\ref{sec:qschur} we introduce quasisymmetric Schur functions, and show in Proposition~\ref{prop:Zbasis} that they form a $\mathbb{Z}$-basis for the algebra of quasisymmetric functions. Section~\ref{sec:qsprops} derives expansions for quasisymmetric Schur functions in terms of monomial and fundamental quasisymmetric functions in Theorems~\ref{the:qsasm} and \ref{the:qsasf}. In Section ~\ref{sub:transmat} we reinterpret these expansions as transition matrices to facilitate the expression of arbitrary quasisymmetric functions in terms of the quasisymmetric Schur function basis. Our main  result of this section, however, is Theorem~\ref{the:pieriqs} in which we give a Pieri rule for quasisymmetric Schur functions. Finally, in Section~\ref{sec:conc} 
we show how to insert the parameter $t$ into our model, defining new quasisymmetric functions which decompose Hall-Littlewood
polynomials; contrast this result with an alternate decomposition obtained by letting $q=0$ in a
formula for Macdonald symmetric functions as a sum of Gessel's fundamental quasisymmetric
functions occurring in \cite{Haglund}, and discuss further avenues to pursue.
% and show the existence of quasisymmetric Hall-Littlewood polynomials.

\section{Quasisymmetric and symmetric functions}\label{sec:qsymsym}

\subsection{Compositions and partitions}\label{sub:comppar}
A \emph{weak composition} $\gamma = (\gamma _1, \gamma _2, \ldots, \gamma _k)$ of $n$, often denoted $\gamma \vDash n$, is a list of nonnegative integers whose sum is $n$. We call the $\gamma _i$ the \emph{parts} of $\gamma$ and $n$ the \emph{size} of $\gamma$, denoted $|\gamma|$. If $\gamma _i$ appears $n_i$ times we abbreviate this subsequence to $\gamma _i^{n_i}$. The \emph{foundation} of $\gamma$ is the set
$$\fo(\gamma)=\{ i\ |\ \gamma _i >0\}.$$If every part of $\gamma$ is positive then we call $\gamma$ a \emph{composition} and call $k:=\ell(\gamma)$ the \emph{length} of $\gamma$. Observe that every weak composition collapses to a composition $\alpha(\gamma)$, which is obtained by removing all $\gamma _i = 0$ from $\gamma$. If every part of $\gamma$ is positive and satisfies $\gamma _1\geq \gamma _2\ge \cdots \geq\gamma _k$ we call $\gamma$ a \emph{partition} of $n$, denoted $\gamma \vdash n$. Observe that every weak composition $\gamma$ determines a partition $\lambda (\gamma)$, which is obtained by reordering the positive parts of $\gamma$ in weakly decreasing order.

\begin{example}
$$\gamma = (3,2,0,4,2,0),\ \fo=\{ 1,2,4,5\},\ \alpha(\gamma)=(3,2,4,2),\ \lambda(\gamma)=(4,3,2,2).$$
\end{example}

Restricting our attention to compositions, there exist three partial orders in which we will be interested. First, given compositions $\alpha, \beta$ we say that $\alpha$ is a \emph{coarsening} of $\beta$ (or $\beta$ is a \emph{refinement} of $\alpha$), denoted $\alpha \succeq \beta$, if we can obtain $\alpha$ by adding together adjacent parts of $\beta$, for example, $(3,2,4,2) \succeq (3,1,1,1,2,1,2)$. Second, we say that $\alpha$ is \emph{lexicographically greater} than $\beta$, denoted $\alpha > _{lex} \beta$,  if $\alpha = (\alpha _1, \alpha _2, \ldots )\neq (\beta _1, \beta _2, \ldots )=\beta$ and the first $i$ for which $\alpha _i \neq \beta _i$ satisfies $\alpha _i > \beta _i$. Third, we say $\alpha \btr \beta$ if $\lambda (\alpha) > _{lex} \lambda (\beta)$ or  $\lambda (\alpha)=\lambda (\beta)$ and $\alpha > _{lex} \beta$.  For example, when $n=4$ we have $$(4)\btr (3,1)\btr (1,3)\btr (2,2)\btr (2,1,1)\btr (1,2,1) \btr (1,1,2) \btr (1,1,1,1).$$

Additionally, to any composition $\beta = (\beta _1, \ldots, \beta _k)$ there is another closely related composition $\beta ^\ast = (\beta _k, \ldots, \beta _1)$, called the \emph{reversal} of $\beta$. Lastly, any composition $\beta = (\beta _1, \beta _2, \ldots, \beta _k)\vDash n$ corresponds to a subset $S(\beta)\subseteq [n-1]=\{1, \ldots , n-1\}$ where
$$S(\beta) = \{\beta _1, \beta _1 +\beta _2, \ldots , \beta _1 +\beta _2 +\cdots + \beta _{k-1}\}.$$Similarly, any subset $S=\{ i_1, i_2, \ldots , i_{k-1}\}\subseteq [n-1]$ corresponds to  a composition $\beta(S)\vDash n$ where
$$\beta(S)= (i_1, i_2 - i_1, i_3 - i_2, \ldots, n-i_{k-1}).$$

\subsection{Quasisymmetric and symmetric function preliminaries}\label{subsec:qsymprelim}
A \emph{quasisymmetric} function is a bounded degree formal power series $F\in \mathbb{Q}[[x_1, x_2, \ldots]]$ such that for all $k$ and $i_1<i_2<\cdots <i_k$ the coefficient of $x_{i_1} ^{\alpha _1}x_{i_2} ^{\alpha _2}\cdots x_{i_k} ^{\alpha _k}$ is equal to the coefficient of $x_{1} ^{\alpha _1}x_{2} ^{\alpha _2}\cdots x_{k} ^{\alpha _k}$ for all compositions $(\alpha _1, \alpha _2, \ldots , \alpha _k)$. The set of all quasisymmetric functions forms a graded algebra $\mathcal{Q}=\mathcal{Q}_0\oplus \mathcal{Q}_1\cdots$.

Two natural bases for quasisymmetric functions are the monomial basis $\{M_{\alpha} \}$ and the fundamental basis $\{F_{\alpha} \}$ indexed by compositions $\alpha = (\alpha_1, \alpha_2, \ldots , \alpha_k)$.    The \emph{monomial} basis consists of $M_0=1$ and all formal power series 
\begin{displaymath}
M_{\alpha} = \sum_{i_1 < i_2 < \cdots < i_k} x_{i_1}^{\alpha_1} x_{i_2}^{\alpha_2} \cdots x_{i_k}^{\alpha_k}.
\end{displaymath}  The \emph{fundamental} basis consists of $F_0=1$ and all formal power series
\begin{displaymath} F_{\alpha} = \sum_{\alpha \succeq \beta} M_{\beta}.\end{displaymath} 
Furthermore, $\mathcal{Q}_n = \rm{span}_{\mathbb{Q}} \{ M_{\alpha} | \alpha \vDash n \} = \rm{span}_{\mathbb{Q}} \{ F_{\alpha} | \alpha \vDash n \}$. We define the algebra of symmetric functions $\Lambda = \Lambda _0 \oplus \Lambda _1 \cdots$ as the subalgebra of $\mathcal{Q}$ spanned by the {\it monomial symmetric functions} $m_0=1$ and all formal power series
\begin{displaymath} m_{\lambda} = \sum_{\alpha\, :\, \lambda(\alpha)=\lambda} M_{\alpha}, \; \; \; \lambda \vdash n  > 0.\end{displaymath} Moreover, we have $\Lambda_n = \Lambda \cap \mathcal{Q}_n$.

\begin{example}
$$F_{(1,2)}=M_{(1,2)}+M_{(1,1,1)}, \quad m_{(2,1)}=M_{(2,1)}+M_{(1,2)}.$$
\end{example}

Perhaps the most well known basis for $\Lambda$ is the basis of Schur functions, $\{ s_{\lambda}\}$, whose definition we devote the next section to.

\section{Schur functions}\label{sec:schur}

\subsection{Diagrams and reversetableaux}\label{sub:diags}
Given a partition $\lambda =(\lambda _1, \lambda _2, \ldots , \lambda _k)$, its corresponding (\emph{Ferrers}) \emph{diagram} is the array of left justified boxes or \emph{cells} with $\lambda _i$ cells in the $i$-th row \emph{from the top}. We abuse notation by using $\lambda $ to refer to both the partition $\lambda$ and its corresponding diagram. We also describe cells by their row and column indices. Given two diagrams $\lambda, \mu$, we say $\mu \subseteq \lambda$ if $\mu _i \leq \lambda _i$ for all $1\leq i \leq \ell(\mu)$, and if $\mu \subseteq \lambda$ then the \emph{skew diagram} $\lambda /\mu$ is the array of cells contained in $\lambda$ but not contained in $\mu$.  In terms of row and column indices
$$\lambda /\mu = \{ (i,j)\ | \ (i,j)\in \lambda, (i,j)\not\in \mu\}.$$The number of cells in $\lambda /\mu$ is called the \emph{size} and is denoted $|\lambda /\mu|$. Two types of skew diagram that will be of particular interest to us later are horizontal strips and vertical strips. We say a skew diagram is a \emph{horizontal strip} if no two cells lie in the same column, and is a \emph{vertical strip} if no two cells lie in the same row.

\begin{example} If
$$\lambda = \tableau{\ &\ &\ &\ \\\ &\ &\ \\\ &\ \\\ &\ }\ , \ \mu = \tableau{\ &\ &\ \\\ &\ \\\ &\ }\ , \ \rho = \tableau{
\ &\ &\ &\ \\\ &\ \\\ \\\ }$$then $\lambda/\mu$ is a horizontal strip and $\lambda / \rho$ is a vertical strip:
$$\lambda / \mu = \tableau{&&&\ \\ &&\ \\ \\\ &\ }\ , \ \lambda / \rho = \tableau{ &\ \\\ \\\ }\ .$$
\end{example}

Reversetableaux are formed from skew diagrams in the following way. Given a skew diagram $\lambda / \mu$ we define a \emph{reversetableau} (or \emph{reverse semistandard Young tableau}), $T$, of \emph{shape} $\lambda / \mu$ to be a filling of the cells with positive integers such that
\begin{enumerate}
\item the entries in the rows of $T$ weakly decrease when read from left to right,
\item the entries in the columns of $T$ strictly decrease when read from top to bottom.
\end{enumerate}
If $|\lambda /\mu|=n$ and the entries are such that each of $1, \ldots , n$ appears once and only once, then we call $T$ a \emph{standard reversetableau}. Classically, given a standard reversetableau, $T$, its \emph{descent set} $D(T)$ is the set of all $i$ such that $i+1$ appears in a higher row. However, by the definition of reversetableau it follows that $i+1$ can only appear
\begin{itemize}
\item strictly above and weakly right
\item weakly below and strictly left
\end{itemize}of $i$. Hence 
$$D(T) = \mbox{the set of all $i$ that do \emph{not} have $i+1$ appearing strictly left of $i$.}$$

\begin{example} In the following reversetableau, to compute $D(T)$ note that 3 is not strictly left of 2.
$$T=\tableau{3&1\\2}\quad ,\quad\ D(T)=\{2\}.$$
\end{example}

The \emph{weight} of a reversetableau, $T$, is the weak composition $w(T)=( w_1(T), w_2(T), \ldots )$ where $w_i(T)=$ the number of times $i$ appears in $T$. The monomial associated with a reversetableau, $T$, is
$$x^T=x_1^{w_1(T)}x_2^{w_2(T)}\cdots .$$For example, the monomial associated with any standard reversetableau, $T$, with $n$ cells is $x^T=x_1x_2\cdots x_n$. We are now ready to define Schur functions.

\subsection{Schur function preliminaries}\label{sub:schur}
There are many ways to define Schur functions, and we begin by defining them as generating functions for reversetableaux. For further details we refer the interested reader to \cite[Chapter 7]{ECII}. Let $\lambda$ be a partition. Then the \emph{Schur function} $s_\lambda$ is
$$s_\lambda = \sum _T x^T$$where the sum is over all reversetableaux, $T$, of shape $\lambda$. We now recall two further classical descriptions, which we include in order to compare with their quasisymmetric counterparts later. The first  describes Schur functions in terms of monomial symmetric functions.

\begin{proposition}\label{prop:sasm}
Let $\lambda , \mu$ be partitions. Then the Schur function $s_\lambda$ is
$$s_\lambda = \sum _\mu K _{\lambda \mu} m_\mu$$where $K_{\lambda \mu}=$ the number of reversetableaux of shape $\lambda$ and weight $\mu ^\ast$.
\end{proposition}

The second description is in terms of fundamental quasisymmetric functions.

\begin{proposition}\label{prop:sasf}
Let $\lambda$ be a partition. Then the Schur function $s_\lambda$ is
$$s_\lambda = \sum _T F_{\beta(D(T))}$$where the sum is over all standard reversetableaux, $T$, of shape $\lambda$. Equivalently,
$$s_\lambda = \sum _\beta d_{\lambda \beta} F_\beta$$where $d_{\lambda \beta}=$ the number of standard reversetableaux, $T$,  of shape $\lambda$ such that $\beta(D(T))=\beta$.
\end{proposition}

\begin{example} We compute
\begin{eqnarray*}s_{(2,1)}&=&m_{(2,1)}+2m_{(1,1,1)}\\
&=&F_{(2,1)}+F_{(1,2)}
\end{eqnarray*}from the reversetableaux
$$\tableau{2&2\\1}\ ,\ \tableau{3&2\\1}\ , \ \tableau{3&1\\2}\ .$$
\end{example}

To close this section we recall two classical products of Schur functions, collectively known as the Pieri rule, which we will later refine to a quasisymmetric setting.

\begin{proposition}[Pieri rule for  Schur functions] \label{prop:pieris}
Let $\lambda$ be a partition.  Then
\[  s_{(n)} s_\lambda = \sum_{\mu} s_\mu \]
where the sum is taken over all partitions $\mu$  such that 
\begin{enumerate}
%\item $\beta \vDash \lvert \alpha \rvert + n$,
\item $\delta = \mu /\lambda$ is a horizontal strip,
\item $\lvert \delta \rvert = n$.
%\item $row_{S(\delta)}(\beta)=\alpha$.
\end{enumerate}
Also,
\[s_{(1^n)} s_\lambda   = \sum_{\mu} s_\mu \]
where the sum is taken over all partitions $\mu$ such that 
\begin{enumerate}
%\item $\beta \vDash \lvert \alpha \rvert + n$,
\item $\epsilon = \mu /\lambda$ is a vertical strip,
\item $\lvert \epsilon \rvert = n$.
%\item $col_{M(\epsilon)}(\beta)=\alpha$.
\end{enumerate}
\end{proposition}

\section{Demazure atoms}\label{sec:atoms}

\subsection{Compositions and diagrams}\label{sub:compdiag}
In this section we define an analogue of reversetableaux that arise naturally in the theory of nonsymmetric Macdonald polynomials. Let $\gamma = (\gamma _1, \gamma _2, \dots , \gamma _n)$ be a weak composition. Then its corresponding \emph{augmented diagram},
 $\augdg (\gamma)$, is the array of left justified cells with  $\gamma _i +1$ cells in the $i$-th row from the top.  Furthermore, the cells of the leftmost column are filled with the integers $1, \ldots , n$ in increasing order from top to bottom, and this $0$-th column is referred to as the \emph{basement}.
 
 \begin{example} 
 $$\augdg (1,0,2)= \tableau{1&\ \\2\\3&\ &\ }\ .$$
 \end{example}
Again we refer to cells by their row and column indices, with the basement being column 0. As with diagrams and reversetableaux we fill the remaining cells of an augmented diagram subject to certain conditions and create semistandard augmented fillings. 

Given an augmented diagram $\augdg (\gamma)$, an \emph{augmented filling}, $\sigma$, is an assignment of positive integer entries to the unfilled cells of $\augdg (\gamma)$. A pair of cells $a = (i,j)$ and $b = (i',j')$ are \emph{attacking} if either $j = j'$ or ($j = j'+1$ and $i > i'$).
An augmented filling $\sigma$ is \emph{non-attacking} if $\sigma(a) \neq \sigma(b)$ whenever $a$ and $b$ are attacking cells.

Then three cells $\{a,b,c\} \in \augdg(\gamma)$ are called a {\it type $A$ triple} if
they are situated as follows
$$\tableau{c & a& \\ \\ &b }$$
where $a$ and $b$ are in the same column, possibly with
cells between them, $c$ is immediately left of $a$, and the length of the row containing $a$ and $c$ is
greater than or equal to the length  of the row containing $b$.  We say that the cells $a,b,c$ form a \emph{type $A$ inversion triple} if their entries, ordered from smallest to largest, form a counter-clockwise orientation.  If two entries are equal, then the entry which appears first when the entries are read top to bottom, right to left,  is considered smallest.

Similarly, three cells $\{a,b,c\} \in \augdg(\gamma)$ are a {\it type $B$ triple} if
they are situated as shown
$$\tableau{ a& \\ \\ b & c}$$
where $a$ and $b$ are in the same column, possibly the basement or with cells between them, $c$ is immediately right  of $b$, and the length of the row
containing $b$ and $c$ is strictly greater than the length of the row containing $a$.  We say that the cells $a,b,c$ form a \emph{type $B$ inversion triple} if their entries, when ordered from smallest to largest, form a clockwise orientation.  Again, if two entries are equal, then the entry which appears first when the entries are read top to bottom, right to left, is considered smallest.

Define a {\it semistandard augmented filling (SSAF)} of \emph{shape $\gamma$}  to be a non-attacking augmented filling of $\augdg (\gamma)$ 
such that the entries in each row  are weakly decreasing when read from left to right (termed \emph{no descents}), and every triple is an inversion triple of type $A$ or $B$.

\begin{remark} Note that in \cite{Mason1} it was shown that the triple and no descent conditions guarantee the augmented filling will be non-attacking. However, we include the extra condition for use in later proofs.
\end{remark}

The weight of a SSAF, $F$, is the weak composition  $w(F)=( w_1(F), w_2(F), \ldots )$ where $w_i(F)=$ (the number of times $i$ appears in $F$) $-1 = $ the number of times $i$ appears in $F$ excluding entries in the basement. Again, the monomial associated with a SSAF, $F$, is
$$x^F=x_1^{w_1(F)}x_2^{w_2(F)}\cdots .$$

\begin{example}
$$F = \tableau{1&{\bf 1}\\
2\\
3&{\bf 3}&{\bf3}}\ , \quad x^F= x_1x_3^2.$$
\end{example}

A SSAF, $F$, of shape $\gamma$ is a \emph{standard augmented filling (SAF)} if for $|\gamma| = n$ we have $x^F= \prod _{i=1} ^n x_i$, and $F$  has \emph{descent set} 
\begin{align*}\D(F):=& \mbox{ the set of all $i$ that do not have $i+1$ appearing strictly left of $i$}\\
&\mbox{ (excluding entries in the basement). }\end{align*}

Similarly, compositions give rise to composition tableaux. Given a composition $\alpha = ( \alpha _1, \alpha _2 , \ldots , \alpha _k)$, its corresponding \emph{composition diagram}, also denoted $\alpha$, is the array of left justified cells with $\alpha _i$ cells in the $i$-th row from the top, and its cells are described by row and column indices.

\begin{definition}\label{def:ComT}Given a composition diagram $\alpha = (\alpha _1, \alpha _2, \ldots , \alpha _\ell)$ with largest part $m$, we define a \emph{composition tableau (ComT)}, $T$, of \emph{shape} $\alpha$ to be a filling of the cells of $\alpha$ with positive integers such that
\begin{enumerate}
\item the entries in the rows of $T$ weakly decrease when read from left to right,
\item the entries in the leftmost column of $T$ strictly increase when read from top to bottom.
\item \emph{Triple rule:} Supplement $T$ by adding enough cells with zero valued entries to the end of each row so that the resulting supplemented tableau, $\hat{T}$, is of rectangular shape $\ell\times m$. Then for $1\leq i<j\leq \ell, 2\leq k \leq m$
$$\left(   \hat{T}(j,k) \neq 0 \mbox{ and } \hat{T}(j,k)\geq \hat{T}(i,k)\right) \Rightarrow  \hat{T}(j,k)> \hat{T}(i,k-1).$$
\end{enumerate}\end{definition}
In exact analogy with reversetableaux, the
weight of a ComT, $T$, is the weak composition $w(T)=( w_1(T), w_2(T), \ldots )$ where $w_i(T)=$ the number of times $i$ appears in $T$. The monomial associated with a ComT, $T$, is
$$x^T=x_1^{w_1(T)}x_2^{w_2(T)}\cdots .$$ Also, a ComT  with $n$ cells is \emph{standard} if  $x^T= \prod _{i=1} ^n x_i$, and has \emph{descent set} $\D(T):=$ the set of all $i$ that do not have $i+1$ appearing strictly left of $i$.

\begin{example} We use a standard composition tableau (ComT) to illustrate our definitions.
$$T = \tableau{5&4&3&1\\6\\8&7&2}\ ,\ \hat{T} = \tableau{5&4&3&1\\6&0&0&0\\8&7&2&0}\ , \ x^T=x_1x_2x_3x_4x_5x_6x_7x_8\ , \ \D(T) = \{2,5,6\}.$$
\end{example}

It transpires that SSAFs and ComTs are closely related, and this relationship will be vital in simplifying subsequent proofs.

\begin{lemma}\label{lem:SSAFandComT} There exists a natural weight preserving bijection between the set of ComTs of shape $\alpha$ and the set of SSAFs of shape $\gamma$ where $\alpha(\gamma)=\alpha$.
\end{lemma}

\begin{example} The following pair consisting of a ComT and SSAF illustrates the natural bijection between them.
$$\tableau{5&4&3&1\\
6\\
8&7&2}\quad \longleftrightarrow\quad \tableau{
1\\
2\\
3\\
4\\
5&{\bf 5}&{\bf 4}&{\bf 3}&{\bf 1}\\
6&{\bf 6}\\
7&\\
8&{\bf 8}&{\bf 7}&{\bf 2}}.$$
\end{example}

\begin{proof}
The mapping that is claimed to be a bijection is clear: given a SSAF, eliminate the basement and any zero parts. For the inverse mapping, given a ComT, let $c$ be the largest  element in the first column.
Allocate a bare basement with $c$ rows.
Place each row of the original ComT to the immediate right of the basement entry that matches the largest row entry.
We need to show that 
\begin{enumerate}
\item the resulting potential ComT satisfies the three rules above, 
\item taking a ComT and applying the inverse operation results in a SSAF.
\end{enumerate}

For the first direction,
assume that $F$ is a SSAF of shape $\gamma$, and that $\sigma$ is the resulting potential ComT of shape $\mu=\alpha (\gamma)$ with $\ell$ rows and $m$ columns and maximum entry $n$.
We first note that  $\sigma$ satisfies Rule 1.
Showing that Rule 2 is satisfied is equivalent to showing that column 1 of $F$ (the column adjacent to the basement) is strictly increasing top to bottom.
Since $F$ is non-attacking, we have that all the entries in each column (in particular, column 1) are distinct.
Note that an entry in column 1 of the SSAF $F$ having value $i$ resides in the cell $(i,1)$.
This follows immediately
since $F$ has no descents and is non-attacking.
Thus it follows that the entries in column 1 of $F$ are strictly increasing,
and so $\sigma$ satisfies Rule 2.

To show that $\sigma$ satisfies Rule 3, suppose to the contrary that there exists a triple of indices $\hat{i},\hat{j},k$ such that
$1 \leq \hat{i} < \hat{j} \leq \ell$, $2 \leq k \leq m$ such that 
$\widehat{\sigma}(\hat{j},k) \neq 0$, $\widehat{\sigma}(\hat{j},k) \geq \widehat{\sigma}(\hat{i},k)$,
and $\widehat{\sigma}(\hat{j},k) \leq \widehat{\sigma}(\hat{i},k-1)$.
Without loss of generality, we may assume that $k$ is minimal over all such triples of indices.
Let $i,j$ be the rows of $F$ corresponding to the respective rows $\hat{i},\hat{j}$ of $\sigma$.
Note that $i < j$.
We consider two cases.

\

\emph{Case: $\gamma _i \geq \gamma _j$.}
In this case, the cell $(i,k)$ of $F$ is nonempty,
and by supposition $F(i,k) < F(j,k) \leq F(i,k-1)$.
But then the cells $(i,k),(j,k),$ $(i,k-1)$ form a non-inversion type A triple,
contradicting the given that $F$ is a SSAF.

\

\emph{Case: $\gamma _i < \gamma _j$.}
Since $F$ has no descents, $F(j,k) \leq F(j,k-1)$, 
and by supposition $F(j,k) \leq F(i,k-1)$.
The cells $(i,k-1),(j,k-1),(j,k)$ form a type B triple, which must be an inversion triple,
so it follows that $F(i,k-1) > F(j,k-1)$.
Since $i < j$, and since the first column of $F$ is strictly increasing, $F(i,1) < F(j,1)$.
Hence there must exist some $k'$, $1 \leq k' < k-1$ such that
$F(i,k') < F(j,k')$ and $F(i,k'+1) > F(j,k'+1)$.
Since $F$ has no descents, $F(i,k'+1) \leq F(i,k')$ and $F(j,k'+1) \leq F(j,k')$, 
which also implies $F(j,k'+1) \leq F(i,k')$.
However, then we see that the cells $(i,k'),(j,k'),(j,k'+1)$ form a non-inversion type B triple,
contradicting the given that $F$ is a SSAF.

\

Thus in both cases we have a contradiction.
It follows that there is no such triple of indices $\hat{i},\hat{j},k$, 
hence $\sigma$ satisfies Rule 3 as well as Rules 1 and 2,
and hence is a ComT.

For the second direction,
assume that $\sigma$ is a ComT, say of shape $\mu$, 
and let $F$ be  obtained by the inverse mapping described above,
which must necessarily be of some shape  $\gamma$, with $\mu = \alpha(\gamma)$.
Since $\sigma$ satisfies Rule 1, $F$ has no descents.
Since $\sigma$ satisfies Rule 2, the first column of $F$ is strictly increasing top to bottom,
and in fact by construction, if cell $(i,1)$ of $F$ is not empty, then $F(i,1) = i$.
In conjunction with this, since $\sigma$ satisfies Rule 3, we have that $F$ must be non-attacking.

Suppose the cells $(i,k),(j,k),(i,k-1)$, $i < j$ form a type A triple in $F$.
If $k = 1$, then $F(i,k) = F(i,k-1) = i < j = F(j,k)$, and so the triple is an inversion triple.
Otherwise $k \geq 2$, and since $\sigma$ satisfies Rules 3 and 1, 
we have that either $F(j,k) < F(i,k) \leq F(i,k-1)$ or $F(i,k) \leq F(i,k-1) < F(j,k)$,
and in both cases the triple is an inversion triple.
Thus all type A triples are inversion triples.

Suppose the cells $(i,k),(j,k),(j,k+1)$, $i < j$ form a type B triple in $F$.
Then $\gamma _i < \gamma _j$.
If $k = 0$, then $F(i,k) = i < j = F(j,k+1) = F(j,k)$, and so the triple is an inversion triple.
Otherwise $k \geq 1$.  
Suppose the triple is not an inversion triple.
This can only happen if  $F(j,k+1) \leq F(i,k) < F(j,k)$.
Let $\hat{i},\hat{j}$ be the rows of $\sigma$ corresponding respectively to the rows $i,j$ of $F$.
Then  $\widehat{\sigma}(\hat{j},k+1) \leq \widehat{\sigma}(\hat{i},k)$,
and Rule 3 then implies that $\widehat{\sigma}(\hat{i},k+1) > \widehat{\sigma}(\hat{j},k+1)$.
Since $\gamma _i < \gamma _j$,
we have $\widehat{\sigma}(\hat{i},\gamma _j) = 0 < \widehat{\sigma}(\hat{j},\gamma _j)$.
There must then exist some $k'$, $k+1 \leq k' < \gamma _j$ such that 
$\widehat{\sigma}(\hat{i},k') > \widehat{\sigma}(\hat{j},k')$ and
$\widehat{\sigma}(\hat{i},k'+1) < \widehat{\sigma}(\hat{j},k'+1)$.
But then we have $\widehat{\sigma}(\hat{i},k'+1) < \widehat{\sigma}(\hat{j},k'+1) \leq \widehat{\sigma}(\hat{j},k') < \widehat{\sigma}(\hat{i},k')$, violating Rule 3.
Thus the triple must be an inversion triple,
and so all type B triples are inversion triples.

We have that $F$ is a non-attacking augmented filling with no descents and in which all type A and type B triples are inversion triples, i.e. $F$ is a SSAF.
\end{proof}

\subsection{Demazure atom preliminaries}\label{sub:dem}
Demazure atoms are formal power series $F\in \mathbb{Q}[[ x_1, x_2, \ldots ]]$, which can be defined combinatorially as follows.

\begin{definition}
Let $\gamma$ be a weak composition. Then the \emph{Demazure atom}, $\da _\gamma$ is
$$\da _\gamma = \sum _F x^F$$where the sum is over all SSAFs, $F$, of shape $\gamma$.
Equivalently, 
$$\da _\gamma = \sum _F x^F$$where the sum is over all ComTs, $F$, of shape $\alpha(\gamma)$ and first column entries $\fo (\gamma)$.
\end{definition}

Note the second definition follows immediately from Lemma~\ref{lem:SSAFandComT}. 

\begin{example}
We compute
$$\da _{(1,0,2)} = x_1x_2x_3 +x_1x_3^2 $$from the SSAFs
$$\tableau{1&{\bf 1}\\2\\3&{\bf 3}&{\bf 2}}\ ,\ \tableau{1&{\bf 1}\\2\\3&{\bf 3}&{\bf 3}}$$or, equivalently, the ComTs
$$\tableau{1\\ 3&2}\ ,\ \tableau{1\\3&3}\ .$$
\end{example}

It transpires that Demazure atoms can be used to describe Schur functions \cite{LS, Mason1}.

\begin{proposition}\label{prop:daschur}
Let $\lambda$ be a partition. Then the Schur function is
$$s_\lambda = \sum _{\gamma\, :\, \lambda (\gamma)=\lambda} \da _\gamma$$where the sum is over all weak compositions $\gamma$.
\end{proposition}

\subsection{Bijection between reversetableaux and SSAFs}
We conclude this section by recalling the bijection $\rho ^{-1}$ from reversetableaux to SSAFs \cite{Mason1}, which we describe algorithmically.

\

Given a reversetableau, $T$, we create a $SSAF$, $\rho ^{-1} (T)=F$, as follows.
\begin{enumerate}
\item If the maximum entry in $T$ is $n$ then allocate a basement with $n$ rows.
\item Taking the entries in $T$ in the first column from top to bottom, place them in  column $k=1$ to the right of the basement  in the uppermost or highest row $i$ of $F$ in which cell $(i,k)$ of $F$ is empty (that is, not yet filled from some earlier column entry of $T$)%the resulting row is weakly decreasing from left to right.
\begin{itemize}
 \item such that the cell $(i,k-1)$ to the immediate left is filled, possibly a basement cell if $k=1$, and 
 \item such that the placement results in no descent. 
 \end{itemize}
\item Repeat with the entries in $T$ in the  column $k$, from top to bottom, placing them in the  column $k$ to the right of the basement for $k = 2, 3, \ldots$.
\end{enumerate}

Eliminating the basement and zero parts from $\rho ^{-1} (T)$ yields a bijection between reversetableaux and ComTs, which we also refer to as $\rho ^{-1}$.

\begin{example}
If $T = \tableau{8&7&3&1\\6&4&2\\5}$ then
$$\rho ^{-1}(T)\quad =\quad \tableau{
1\\
2\\
3\\
4\\
5&{\bf 5}&{\bf 4}&{\bf 3}&{\bf 1}\\
6&{\bf 6}\\
7&\\
8&{\bf 8}&{\bf 7}&{\bf 2}} \quad \equiv\quad  \tableau{5&4&3&1\\
6\\
8&7&2}\ .$$
\end{example}

\section{Quasisymmetric Schur functions}\label{sec:qschur}
We now define our main objects of study and derive some elementary properties about them.

\begin{definition}\label{def:qschur}
Let $\alpha$ be a composition. Then the \emph{quasisymmetric Schur function} is
$$\qs _\alpha = \sum _{\gamma\, :\, \alpha (\gamma)=\alpha} \da _\gamma$$where the sum is over all weak compositions $\gamma$.
\end{definition}

\begin{example} Restricting ourselves to three variables we compute
\begin{eqnarray*}\qs _{(1,2)} &=& \da _{(1,2,0)}+\da _{(1,0,2)} + \da _{(0,1,2)}\\
&=& x_1x_2^2 + x_1x_2x_3+ x_1x_3^2  + x_2x_3^2\end{eqnarray*}where the summands arise from all SSAFs of shape $(1,2,0), (1,0,2)$ and $(0,1,2)$:
\begin{displaymath}
\tableau{1&{\bf 1}\\ 2&{\bf 2}&{\bf 2}\\3} \; \; \; \; \;  \; \; \; \; \; \; \; \tableau{ 1&{\bf 1}\\2\\ 3&{\bf 3}&{\bf 2}} \; \;  \; \; \; \tableau{ 1&{\bf 1}\\2\\ 3&{\bf 3}&{\bf 3}} \; \; \; \; \; \; \; \; \; \; \; \; \tableau{ 1\\2&{\bf 2}\\ 3&{\bf 3}&{\bf 3}}
\end{displaymath}or, equivalently, from ComTs
\begin{displaymath}
\tableau{1\\2&2} \; \; \; \; \;  \; \; \; \; \; \; \; \tableau{ 1\\3&2} \; \;  \; \; \; \tableau{ 1\\3&3} \; \; \; \; \; \; \; \; \; \; \; \; \tableau{ 2\\3&3}\ .
\end{displaymath}
\end{example}

As illustrated by this example, we shall see later that the functions are indeed quasisymmetric, but first we focus on their connection to Schur functions. 

Recall from Proposition~\ref{prop:daschur} that the Schur function $s_\lambda$ decomposes into the sum of all $\da _\gamma$ such that $\lambda (\gamma)=\lambda$. Hence by Definition~\ref{def:qschur} we obtain the decomposition of the Schur function in terms of quasisymmetric functions
$$s_\lambda = \sum _{\alpha\, :\, \lambda (\alpha)=\lambda} \qs _\alpha ,$$which immediately evokes the definition of \emph{monomial} symmetric functions in terms of \emph{monomial} quasisymmetric functions
$$m_\lambda = \sum _{\alpha\, :\, \lambda (\alpha)=\lambda} M _\alpha .$$Thus, the parallel construction justifies the use of the word Schur. We also prove the functions are quasisymmetric by describing quasisymmetric Schur functions in terms of fundamental quasisymmetric functions.

\begin{proposition}\label{prop:qsasf} 
Let $\alpha$ be a composition. Then 
$$\qs_\alpha = \sum _T F_{\beta(D(T))}$$where the sum is over all standard reversetableaux, $T$, of shape $\lambda(\alpha)$ that map under $\rho ^{-1}$ to a SSAF of shape $\gamma$ satisfying $\alpha(\gamma)=\alpha$ (or, equivalently, under $\rho ^{-1}$ to a ComT of shape $\alpha$).
\end{proposition}

\begin{proof}
To prove this we need to show
\begin{enumerate}
\item for each $T$ satisfying the conditions stated, $F_{\beta(D(T))}$ is a summand of $\qs _\alpha$ appearing exactly once,
\item the coefficient of each monomial appearing in $\qs _\alpha$ is equal to the sum of its coefficients in each of the $F_{\beta(D(T))}$s in which it appears.
\end{enumerate}
To show the first point note that $F_{\beta(D(T))}$ is a sum of monomials that arise from reversetableaux, which standardize to $T$. Furthermore, any reversetableau $\tilde{T}$ that standardizes to $T$, denoted $std(\tilde{T})=T$, maps under $\rho ^{-1}$ to a SSAF 
$\rho ^{-1} (\tilde{T})$ that standardizes to $\rho ^{-1} ({T})$ \cite{Mason1}. That is, if we say for a SSAF, $F$, that its standardization is $\rho^{-1}(std(\rho(F)))$ then $\rho ^{-1}(T)= \rho ^{-1}(std(\rho(\rho ^{-1} (\tilde{T})))$, where $\rho$ is the inverse of $\rho ^{-1}$. Thus, if $T$ is of shape $\lambda (\alpha)$ and $\rho ^{-1} ({T})$ is of shape $\gamma$ such that $\alpha(\gamma)=\alpha$, then $F_{\beta(D(T))}$ is a summand of $\qs _\alpha$ appearing exactly once.

To show the second point, observe if given a SSAF $\rho ^{-1} (\tilde{T})$ of shape $\gamma$ such that $\alpha (\gamma)=\alpha$, which contributes a monomial towards $\qs _\alpha$ and also standardizes to $\rho ^{-1} ({T})$, then under  $\rho $ this maps bijectively  to $\tilde{T}$ that standardizes to $T$ of shape $\lambda(\alpha)$. Computing $D(T)$ then yields which fundamental quasisymmetric function the monomial belongs to.
\end{proof}

A combinatorially more straightforward description in terms of the $F_\alpha$ is given in the next section, and hence we delay giving an example until then. We will now show that, in fact, the set of all quasisymmetric Schur functions forms a basis for $\mathcal{Q}$. Before we do this, we will work towards two lemmas.

For a composition $\alpha$, let $T_\alpha$ be the unique standard reversetableau of shape $\lambda (\alpha)$ and $\beta(D(T))=\alpha$.  To see that $T_\alpha$ exists, construct the left justified array of cells of shape $\alpha ^\ast$,   which has $1, \ldots, \alpha_1$ in the bottom row and 
$$\alpha_1+\cdots + \alpha _{i-1}+1, \ldots , \alpha_1+\cdots + \alpha _{i}$$ in the $i$-th row from bottom appearing in decreasing order when read from left to right.
Then move every cell as far north as possible to form $T_\alpha$. To see that $T_\alpha $ is unique, note that the number of descents in $T_\alpha$ is one less than the number of rows in $T_\alpha$ and so all entries in the first column except $n$ must be all $i$ such that $i\in D(T_\alpha)$. This and the fact that $T_\alpha$ must be a reversetableau yield uniqueness.

\begin{example}
If $\alpha=(1,3,2)$ then we construct
$  \tableau{ 6&5\\4&3&2\\1}$ and 
$T_\alpha = \tableau{
6&5&2\\
4&3\\
1}$\ .
\end{example}

The following lemma is straightforward to verify using the algorithm for $\rho^{-1}$.

\begin{lemma}\label{lem:TSSAF}
For a composition $\alpha =(\alpha _1,\alpha _2,\ldots ,\alpha _k)\vDash n$, performing $\rho ^{-1}$ on $T_\alpha$ yields the SSAF with basement $1, \ldots , n$ and row $\alpha _1$ containing $1, \ldots , \alpha _1$, row $\alpha _1+\alpha _2$ containing $\alpha _1 +1, \ldots , \alpha _1 +\alpha _2$ etc. Equivalently, performing $\rho ^{-1}$ on $T_\alpha$ yields the ComT with row $1$ containing $1, \ldots , \alpha _1$, row $2$ containing $\alpha _1 +1, \ldots , \alpha _1 +\alpha _2$ etc.
\end{lemma}

\begin{lemma}\label{cor:leadHsummand}
$F_\alpha$ will always be a summand of $\mathcal{S}_\alpha$ with coefficient 1.
\end{lemma}

\begin{proof}
This follows immediately from the existence and uniqueness of $T_{\alpha}$, Proposition~\ref{prop:qsasf} and Lemma~\ref{lem:TSSAF}.
\end{proof}

We are now ready to prove that quasisymmetric Schur functions form a basis for $\mathcal{Q}$.

\begin{proposition}{\label{prop:Zbasis}}
The set $\{ \mathcal{S}_{\alpha} | \alpha \vDash n\}$ forms a $\mathbb{Z}$-basis for $\mathcal{Q}$.
\end{proposition}

\begin{proof}
For a fixed $n$ and $\alpha =(\alpha _1,\ldots, \alpha _{\ell(\alpha)}) \vDash n$ consider the summand $F_\delta$ appearing in $\mathcal{S}_\alpha$. By Proposition \ref{prop:qsasf}  it follows that $\lambda(\alpha)\geq _{lex} \lambda (\delta)$ because if not then the first $i$ when $\lambda (\delta )_i >  \lambda (\alpha)_i$ will yield a row in any diagram $\lambda(\alpha)$ that cannot be filled to create a standard reversetableau, $T$, satisfying $D(T)=S(\delta)$. If $\lambda (\alpha)=\lambda (\delta)$ then by Lemma~\ref{lem:TSSAF} and the uniqueness of $T_{\alpha}$  we know the coefficient of $F_\delta$ will be $0$ unless $\alpha = \delta$.

Let $M$ be the matrix whose rows and columns are indexed by $\alpha \vDash n$ ordered by $\btr$ and entry $M_{\alpha\delta}$ is the coefficient of $F_\delta$ in  $\mathcal{S}_\alpha$.  By the above argument and Lemma~\ref{cor:leadHsummand} we have that $M$ is upper unitriangular, and the result follows.
\end{proof}

\section{Properties of quasisymmetric Schur functions}\label{sec:qsprops}
A natural question to ask about quasisymmetric Schur functions is how many properties of Schur functions refine to \emph{quasisymmetric} Schur functions? In this regard there are many avenues to pursue. In this section we provide the expansion of a quasisymmetric Schur function in terms of monomial symmetric functions, and a more explicit expression in terms of fundamental quasisymmetric functions. Our main result of this section, however, is to show that quasisymmetric Schur functions exhibit a Pieri rule that naturally refines the original Pieri rule for Schur functions.

To appreciate these quasisymmetric refinements we invite the reader to compare the classical Schur function properties of Propositions~\ref{prop:sasm}, \ref{prop:sasf} and \ref{prop:pieris} with the quasisymmetric Schur function properties of Theorems~\ref{the:qsasm}, \ref{the:qsasf} and \ref{the:pieriqs}, respectively.

\begin{theorem}\label{the:qsasm}
Let $\alpha, \beta$ be compositions. Then  
$$\qs _\alpha = \sum _\beta K _{\alpha \beta} M_\beta$$where $K_{\alpha\beta}=$ the number of SSAFs of shape $\gamma$ satisfying $\alpha (\gamma) = \alpha$ and weight $\beta$ (or, equivalently, 
$K_{\alpha\beta}=$ the number of ComTs of shape $\alpha$ and weight $\beta$).
\end{theorem}

\begin{proof}
We know
$$\qs _\alpha = \sum _{\gamma\,:\,\alpha(\gamma)=\alpha} \da _{\gamma}= \sum x^F = \sum _\beta c_{\alpha\beta} M_\beta$$where the middle sum is over all SSAFs $F$ of shape $\gamma$ where $\alpha(\gamma)=\alpha$. The leading term of any $M_\beta$ appearing in the last sum is $x_1^{\beta _1}x_2^{\beta _2}\cdots x_\ell^{\beta _\ell}$, and the number of times it will appear is, by the middle equality, the number of SSAFs of shape $\gamma$ where $\alpha(\gamma)=\alpha$ and weight $\beta$. Hence $c_{\alpha \beta}=K_{\alpha \beta}$ and the result follows.
\end{proof}

\begin{theorem}\label{the:qsasf}
Let $\alpha, \beta$ be compositions. Then 
$$\qs _\alpha = \sum _\beta d_{\alpha \beta} F_\beta$$where $d_{\alpha \beta}=$ the number of SAFs $T$ of shape $\gamma$ satisfying $\alpha (\gamma) = \alpha$  and $\beta(\D(T))=\beta$ (or, equivalently,  $d_{\alpha \beta}=$ the number of standard ComTs $T$ of shape $\alpha$ and $\beta(\D(T))=\beta$).
\end{theorem}

\begin{proof}
Since $\qs _\alpha = \sum _T F_{\beta(D(T))}$ where the sum is over all standard reversetableaux, $T$, of shape $\lambda (\alpha)$ that map under $\rho ^{-1}$ to a SSAF of shape $\gamma$ that satisfies $\alpha(\gamma)=\alpha$, and since $\rho ^{-1}$ maps the entries appearing in column $j$ of $T$ to column $j$ of $\rho ^{-1}(T)$, the result follows.
\end{proof}

\begin{example} We compute
\begin{eqnarray*}
\qs _{(1,2)}&=&M_{(1,2)}+M_{(1,1,1)}\\
&=& F_{(1,2)}
\end{eqnarray*}from the SSAFs
$$\tableau{1&{\bf 1}\\2&{\bf 2}&{\bf 2}}\quad ,\quad \tableau{1&{\bf 1}\\2\\3&{\bf 3}&{\bf 2}}$$or, equivalently, the ComTs
$$\tableau{1\\2&2}\quad ,\quad \tableau{1\\3&2}$$for the first equality, and just the latter SAF or standard ComT for the second equality.
\end{example}

\begin{remark}
For a composition $\alpha=(\alpha _1, \alpha _2, \ldots , \alpha _{\ell(\alpha)})$ \cite{MacM} defines the \emph{symmetric} function indexed by a composition, known as a \emph{ribbon Schur function}, $r _\alpha$. The relationship between ribbon Schur functions and similarly indexed quasisymmetric Schur functions is straightforward to deduce as follows.

By the Littlewood-Richardson rule, say \cite[Chapter 7]{ECII}, we have
$$r_\alpha = \sum c_{\alpha\lambda} s_\lambda$$where $c_{\alpha\lambda}$ is the number of Littlewood-Richardson fillings of the connected skew diagram containing no $2\times2$ skew diagram that has $\alpha _1$ cells in the top row, $\alpha _2$ cells in the second row etc. Since $s_\lambda = \sum _{\alpha\, :\, \lambda (\alpha)=\lambda} \qs _\alpha$ it immediately follows that
$$r_\alpha = \sum c_{\alpha\lambda (\beta)} \qs_\beta$$where $c_{\alpha\lambda (\beta)}$ is   as above.

In \cite[Theorem 4.1]{HDL} necessary and sufficient conditions for equality of ribbon Schur functions were determined. Meanwhile, in \cite[Theorem 2.2]{SvW} necessary and sufficient conditions for uniqueness of Littlewood-Richardson fillings  were proved. Combining these results with the above, we conclude that the simple relationship between Schur functions and quasisymmetric Schur functions is only achieved again with $r_{(u, 1^v)}$ or $r_{(1^v, u)}$, that is
$$r_{(u, 1^v)}=r_{(1^v, u)}= \sum _{\lambda (\alpha)=(u, 1^v)} \qs _\alpha.$$
Thus concludes our remark.
\end{remark}

We now come to our Pieri rule for quasisymmetric Schur functions, whose proof we delay until the next subsection, and whose statement requires the following definitions.

\begin{remark} In  practice the following $rem _s$ operator subtracts 1 from the rightmost part of size $s$ in a composition, or returns the empty composition. Meanwhile the $row _{\{ s_1 < \cdots < s_j \}}$ operator subtracts 1 from the rightmost part of size $s_j, s_{j-1}, \ldots$ recursively. Similarly, the $col _{\{m_1 \leq \cdots \leq m_j \}}$ operator subtracts 1 from the rightmost part of size $m_1, m_2, \ldots$ recursively. 
\end{remark}

\begin{example} If $\alpha = (1,2,3)$ then
$$row _{\{2,3\}}(\alpha)= rem_2(rem_3((1,2,3)))= rem_2((1,2,2))= (1,2,1)$$and
$$col _{\{2,3\}}(\alpha)=rem_3(rem_2((1,2,3)))= rem_3((1,1,3))= (1,1,2).$$
\end{example}

Now we define these three operators formally. Let $\alpha = (\alpha_1,\ldots,\alpha_k)$ be a composition whose largest part is $m$, and let $s \in [m] $.
If there exists  $1\leq i\leq k$ such that $s=\alpha_i$ and $s\neq\alpha_j$ for all $j>i$, then define
\[ rem_s(\alpha)=(\alpha_1,\ldots,\alpha_{i-1},(s-1),\alpha_{i+1},\ldots,\alpha_k), \]
otherwise define $rem_s(\alpha)$ to be the empty composition.
Let $S = \{ s_1 < \cdots < s_j \}$.  Then define
\[ row_S(\alpha)=rem_{s_1}(\ldots(rem_{s_{j-1}}(rem_{s_j}(\alpha)))\ldots). \]
Similarly let $M = \{m_1 \leq \cdots \leq m_j \}$.  Then define
\[ col_M(\alpha)=rem_{m_j}(\ldots(rem_{m_2}(rem_{m_1}(\alpha)))\ldots). \]We collapse $row_S(\alpha)$ or $col_M(\alpha)$ to obtain a composition if needs be.

%Here we will use the notation $\lambda(\alpha)$ in place of the notation $sort(\alpha)$ used in an earlier section to denote the partition obtained from the composition $\alpha$ by sorting the parts of $\alpha$ into weakly decreasing order.
For any horizontal strip $\delta$ we denote by $S(\delta)$ the set of columns its skew diagram occupies,
and for any vertical strip $\epsilon$ we denote by $M(\epsilon)$ the multiset of columns  its skew diagram occupies, where multiplicities for a column are given by the number of cells in that column. We are now ready to state our refined Pieri rule.

\begin{theorem}[Pieri rule for quasisymmetric Schur functions] \label{the:pieriqs}
Let $\alpha$ be a composition.  Then
\[ \mathcal{S}_{(n)} \mathcal{S}_\alpha  = \sum_{\beta} \mathcal{S}_\beta  \]
where the sum is taken over all compositions $\beta$  such that 
\begin{enumerate}
%\item $\beta \vDash \lvert \alpha \rvert + n$,
\item $\delta = \lambda(\beta)/\lambda(\alpha)$ is a horizontal strip,
\item $\lvert \delta \rvert = n$,
\item $row_{S(\delta)}(\beta)=\alpha$.
\end{enumerate}
Also,
\[\mathcal{S}_{(1^n)} \mathcal{S}_\alpha   = \sum_{\beta} \mathcal{S}_\beta \]
where the sum is taken over all compositions $\beta$ such that 
\begin{enumerate}
%\item $\beta \vDash \lvert \alpha \rvert + n$,
\item $\epsilon = \lambda(\beta)/\lambda(\alpha)$ is a vertical strip,
\item $\lvert \epsilon \rvert = n$,
\item $col_{M(\epsilon)}(\beta)=\alpha$.
\end{enumerate}
\end{theorem}

For a more visual interpretation of Theorem~\ref{the:pieriqs} we use composition diagrams in place of compositions in the next example. Then $rem _s$ is the operation that removes the rightmost cell from the lowest row of length $s$.

\begin{example}If we place $\bullet$ in the cell to be removed then
$$rem_1((1,1,3))= 
\tableau{\ \\
\bullet\\
\ &\ &\ &}= (1,3) .$$If we wish to compute $\mathcal{S} _{(1)}\mathcal{S} _{(1,3)}$ then we consider the four skew diagrams 
$$(4,1)/(3,1),\ (3,2)/(3,1),\ (3,1,1)/(3,1),\ (3,1,1)/(3,1)\ (again)$$with horizontal strips containing one cell in column $4,2,1,1$ respectively. Then
$$row _{\{4\}} ((1,4)) = \tableau{\ \\
\ &\ &\ &\bullet} \quad
row _{\{2\}} ((2,3)) = \tableau{\ & \bullet \\
\ &\ &\  }$$
$$row _{\{1\}} ((1,3,1)) = \tableau{\ \\
\ &\ &\  \\
\bullet} \quad 
row _{\{1\}} ((1,1,3)) = \tableau{\ \\
\bullet\\
\ &\ &\ } \ $$and hence
$$\mathcal{S}_{(1)}\mathcal{S}_{(1,3)}=
\mathcal{S}_{(1,4)} + \mathcal{S}_{(2,3)} + \mathcal{S}_{(1,3,1)} +\mathcal{S}_{(1,1,3)}.$$
\end{example}

Classically, the Pieri rule gives rise to Young's lattice on partitions in the following way. Let $\lambda, \mu$ be partitions, then $\lambda$ covers $\mu$ in Young's lattice if the coefficient of $s_\lambda$ in $s_{(1)}s_\mu$ is 1. Therefore, Theorem~\ref{the:pieriqs} analogously gives rise to a poset on compositions: Let $\alpha, \beta$ be compositions, then $\beta$ covers $\alpha$ if the coefficient of $\qs _\beta$ in $\qs _{(1)} \qs _\alpha$ is 1. It would be interesting to see what properties of Young's lattice are exhibited by this new poset, which differs from the poset of compositions in \cite{BergBousDul}, and contains Young's lattice as a subposet.

\subsection{Proof of the Pieri rule for quasisymmetric Schur functions}
In order to prove our Pieri rule we require three known combinatorial constructs, which we recall here in terms of reversetableaux for convenience.

The first construct is \emph{Schensted insertion}, which inserts a positive  integer $k_1$ into a reversetableau $T$, denoted $T\leftarrow k_1$ by
\begin{enumerate}
\item if $k_1$ is less than or equal to the last entry in  row 1, place it there, else
\item find the leftmost entry in that row strictly smaller than $k_1$, say $k_2$, then
\item replace $k_2$ by $k_1$, that is $k_1$ \emph{bumps} $k_2$.
\item Repeat the previous steps with $k_2$ and row 2, $k_3$ and row 3, etc.
\end{enumerate}
The set of cells whose values are modified by the insertion, including the final cell added, is called the \emph{insertion path}, and the final cell is called the \emph{new cell}.

\begin{example}
$$\tableau{7 & 5 & 4 & 2  \\
  6 & 4 & 3  \\
    3 & 2 & 2 \\
   1 & 1 \\
}\quad  \gets 5  \qquad
=  \qquad
\tableau{
 7        & 5 & \gx{5} & 2  \\
  6        & 4 & \gx{4}  \\
    3 & \gx{3}  & 2 \\
       \gx{2} & 1 \\
   \gx{1}  \\
}$$
where the bold italic cells indicate the insertion path.
\end{example}

Insertion paths have the useful property encompassed in the next lemma, commonly known as the \emph{row bumping lemma}.

\begin{lemma}[Row bumping lemma] \label{lemma:bumping}
Let $T$ be a reversetableau. Consider two successive insertions $(T\gets x )\gets x'$,
giving rise to two insertion paths $R$ and $R'$, with respective new cells $B$ and $B'$.
\begin{enumerate}
\item If $x \geq x'$, then $R$ is strictly left of $R'$, and $B$ is strictly left of and weakly below $B'$.
\item If $x < x'$, then $R'$ is weakly left of $R$, and $B'$ is weakly left of and strictly  below $B$.
\end{enumerate}
\end{lemma}

The second combinatorial construct we require is the \emph{plactic monoid}, which can be described as the monoid whose elements consist of all reversetableaux. To describe the product, recall the \emph{row reading word} of a reversetableau, $T$, is the sequence of the entries of the cells of $T$ read from left to right, and bottom to top. It is denoted $w_{row}(T)$. For example, $w_{row}(T)= 113226437542$ for the original reversetableau in the previous example. Then, given reversetableaux $T$ and $U$, the plactic monoid product is
$$T\cdot U = ((T \gets w_1) \gets w_2) \cdots \gets w_n$$where $w_{row} (U) = w_1w_2\cdots w_n.$ The empty reversetableau is the monoid identity.

The group ring of the plactic monoid, $R$, is called the \emph{reversetableaux ring} and $S_\lambda \in R$ is
$$S _\lambda = \sum T$$where the sum is over all reversetableaux, $T$, of shape $\lambda$.

There exists a surjective homomorphism
\begin{eqnarray*}\varepsilon\ :\ R&\longrightarrow& \mathbb{Z}[[x_1, x_2, \ldots ]]\\
T&\mapsto & x^T\end{eqnarray*}that importantly satisfies
$$\varepsilon (S_\lambda ) = s_\lambda .$$

The third, and last, construct is an analogy to Schensted insertion for a SSAF, or \emph{skyline insertion}. We state it here for ComTs since ComTs will be used in the remaining proofs. However, it can be found in its original form in \cite[Procedure 3.3]{Mason1}. 

Suppose we start with a ComT $F$ whose longest row has length $r$.
To insert a positive integer $k_1$, the result being denoted $k_1\to
F$,
scan column positions starting with the top position in column $j 
= r+1$.
\begin{enumerate}
\item If the current position is empty and at the end of a row of  
length $j-1$, and $k_1$ is weakly less than
the last entry in the row, then place $k_1$ in this empty position  
and stop. Otherwise, if the position
is nonempty and contains $k_2<k_1$ and $k_1$ is weakly less than the  
entry to the immediate left of $k_2$,
let $k_1$ bump $k_2$, i.e. swap $k_2$ and $k_1$.
\item Using the possibly new $k_i$ value, continue scanning successive  
positions in the column top to bottom as in the previous step, bumping  
whenever possible, and then continue scanning at the top of the next  
column to the left. (Decrement $j$.)
\item If an element is bumped into the first column, then create a new  
row containing one cell to contain the element,
placing the row such that the first column is strictly increasing top
to bottom, and stop.
\end{enumerate}

The set of cells whose values are modified by the insertion, including the final cell added, is called the \emph{insertion sequence}, and the final cell is called the \emph{new cell}. The row in which the new cell is added is called the \emph{row augmented by the insertion}, and we note that the number of cells, or length of the row, increases by one.

\begin{example}
$$5 \to \quad \tableau{
   1 & 1 \\
   3 & 2 & 2 & 2 \\
   6 & 5 & 4  \\
   7 & 4 & 3   \\
}  
 \qquad =  \qquad
\tableau{
   1  & 1\\
   \gx{2}               \\
   3  & \gx{3} & 2 & 2 \\
   6        & 5 & \gx{5}  \\
   7        & 4 & \gx{4}   \\
} $$where the bold italic cells indicate the insertion sequence.
\end{example}

Schensted and skyline insertion commute in the following sense \cite[Proposition 3.1]{Mason1}.

\begin{proposition}\label{prop:commute}
If $\rho$ is the inverse map of $\rho ^{-1}$ and $F$ is a ComT then
$$\rho( k \to F) = (\rho(F) \gets k).$$
\end{proposition}

We are now ready to prove the Pieri rule for quasisymmetric Schur functions after we prove

\begin{lemma} \label{lemma:key}
Let $D$ be a ComT, $k$ a positive integer, and $D' = k \to D$ with row $i$ of $D'$ being the row augmented by the insertion.
Then for all rows $r > i$ of $D'$
 $$\mbox{ length of row }i \neq \mbox{ length of row }r.$$
\end{lemma}

\begin{proof}
Assume that the lemma is false, that is, that there exists a row $r>i$ of $D'$ having the same length as the augmented row $i$, say length equal to $j$.
Note that in this case, the $r$-th row of $D'$ is the same as the $r$-th row of $D$, except in the case that the augmented row $i$ is a new row of length 1, in which case the $(r+1)$-th row of $D'$ is the same as the $r$-th row of $D$.
In the  algorithm  for inserting a new element $k$ into a ComT $D$,
consider the value $x$ that was bumped from column $j+1$ into column $j$.
%that is, $x$ is the value of the variable $k$ in the pseudo-code at the point that the last row of column $j$ has been completely scanned %and $j$ is being decremented.

%On the one hand 
This bumped value $x$ must be larger than $D(r,j)$, for otherwise either $x$ was the value of the variable $k$ compared against $D(r,j)$ during the pass of the algorithm over column (j+1), in which case the value $x$ would have been placed in the vacant position $D(r,j+1)$, 
or $x$ was bumped from position $D(s,j+1)$ for some row $s>r$, in which case $D$ would have violated the triple rule for ComTs (namely $D(r,j)\geq x=D(s,j+1)>D(r,j+1)$), a contradiction in either case.

Now if $j=1$, then $x$ was simply inserted into $D'$ as the new row $i$ of length one. 
However, since the first column is strictly increasing top to bottom, $x$ would have been inserted as a new row \emph{after} (higher row number than)  $D(r,j)$, i.e. $i > r$, contrary to supposition.
So we can assume that $j>1$.

Recall that the entries in any given column are all distinct.
We must have $D'(r,j) >  D'(i,j-1)$, for otherwise we would have had a triple rule violation in $D$ 
(namely $D(i,j-1)\geq D(r,j)>D(i,j)=\text{empty}$).
This then would require that $D'(r,j)>D'(i,j)$ as well.
Now consider the portion of the insertion sequence that lies in column $j$, say in rows $\{i_0,\ldots,i_t=i\}$, whose first value, scanning top to bottom, is $x=D'(i_0,j)$ and whose last value is $D'(i,j)=D'(i_t,j)$.  
Since $x > D(r,j)=D'(r,j) > D'(i,j)$, and since the entries in the insertion path are decreasing top to bottom,  there must be some index $0\leq \ell < t$ such that $D'(i_\ell,j)>D'(r,j) > D'(i_{\ell+1},j)$.
This would imply $D'(i_\ell,j-1)>D'(r,j)$ as well.
However, since $D(i_k,j)=D'(i_{k+1},j)$ for all $0\leq k < t$, and $D(h,k)=D'(h,k)$ for all $k<j$, this would imply a triple rule violation in $D$, namely $D(i_\ell,j-1)>D(r,j) > D(i_\ell,j)$.

Thus in all cases we have a contradiction.
\end{proof}

We  note that  as Schensted insertion for reversetableaux  is reversible (invertible), so the analogous insertion into ComTs is reversible.  
In particular, given a ComT  $D$ of shape $\alpha$ and a given positive integer $\ell = \alpha_i$ for some $i$, where we assume that $i$ is the largest index such that $\ell = \alpha_i$,
then one can \emph{uninsert} an element $k$ from $D$ to obtain a ComT $D'$ such that $D = k \to D'$ and the shape of $D'$ is $(\alpha_1,\ldots,\alpha_{i-1},\ell-1,\alpha_{i+1},\ldots)$, that is, the shape obtained from $D$ by removing the last square from row $i$.

\begin{proof} (of Theorem~\ref{the:pieriqs})
We start with the first formula.
We consider $S_n$ to be the sum of reversetableaux in the reversetableaux ring $R$ of shape $(n)$,
and $H_\alpha$ to be the sum of reversetableaux in $R$ of shape $\lambda(\alpha)$ which map to a ComT of shape $\alpha$ under the mapping $\rho^{-1}$.
We consider a typical term $U\cdot V$ of the product  $H_\alpha \cdot S_n$, where $U$ is one of the reversetableau terms of $H_\alpha$ and $V$ is one of the reversetableau terms of $S_n$.
Suppose the reversetableau $C = U\cdot V$, and ComT $D = \rho^{-1}(C)$,
where the shape of $D$ is $\beta$.
Note that $U$ has shape $\lambda(\alpha)$ and $C$ has shape $\lambda(\beta)$.
The set of $n$ new cells added to $U$ in the product $U\cdot V$ to form $C$ forms the skew reversetableau $C/U$ of shape $\lambda(\beta)/\lambda(\alpha)$.
Now $V$ will be of the form \[ V =\tableau{  x_1 & x_2 & \cdots & x_n} \]
where $x_1 \geq x_2 \geq \cdots \geq x_n$.  
Thus \[ C = U \gets x_1 \gets x_2 \gets \cdots \gets x_n. \]
By Lemma \ref{lemma:bumping}, $C/U$ is a horizontal strip with $n$ cells,
and over the successive  insertions, the cells in this horizontal strip are added to $U$ from left to right,
say in columns $j_1<j_2<\cdots<j_n$.
Suppose ComT $E=\rho^{-1}(U)$, which by assumption has shape $\alpha$.
Under the map $\rho^{-1}$ and insertion for ComTs, 
the corresponding new cells added to $E$ to form $D$ are added to columns $j_1,j_2,\ldots,j_n$ in the same order by Proposition~\ref{prop:commute}.
By Lemma \ref{lemma:key}, each time a new cell is added, the augmented row in which it appears is the last row in the new diagram of that length.
That is, assuming that  $\alpha_i=j_1-1$ and $\alpha_k\neq j_1-1$ for all $k>i$, then the shape of $(E\gets x_1)$ is $\alpha'=(\alpha_1,\ldots,\alpha_{i-1},j_1,\alpha_{i+1},\ldots,\alpha_{\ell(\alpha)})$,
that is $\alpha = rem_{j_1}(\alpha')$.
The pattern continues, that is, if the shape of $(E\gets x_1\gets x_2)$ is $\alpha''$,
then $\alpha' = rem_{j_2}(\alpha'')$, etc.,
and by induction we have 
\[ \alpha = rem_{j_1}(\ldots(rem_{j_{n-1}}(rem_{j_n}(\beta)))\ldots) = row_J(\beta) \]
where $J = \{j_1,\ldots ,j_n\}$.
Thus $C = U\cdot V$ is a term (reversetableau summand) of $H_\beta$ where $\beta$ is one of the summand indices specified by the formula.

%\bigskip
Conversely, suppose the reversetableau $C$ is a term of $H_\beta$ where $\beta$ is one of the summand indices on the right hand side of the formula. 
By definition, $\lambda(\beta)/\lambda(\alpha)$ is a horizontal strip with $n$ cells,
say in columns $j_1<j_2<\cdots<j_n$.
Since  insertion is reversible, we can perform uninsertion on $C$, removing the cells of the horizontal strip starting with the last column $j_n$ and working left.
Uninserting the bottom cell from column $j_n$ yields an element $x_n$, then uninserting the bottom cell from column $j_{n-1}$ yields an element $x_{n-1}$, etc.
\[ (\ldots ((C \xrightarrow{j_n} x_n) \xrightarrow{j_{n-1}} x_{n-1}) \ldots ) \xrightarrow{j_1} x_1 \]
Let $U$ be the reversetableau resulting from uninserting the $n$ cells.
Now Lemma  \ref{lemma:bumping} implies that $x_1 \geq x_2 \geq \cdots \geq x_n$,
and so we may set $V$ to be the reversetableau of shape $(n)$ having these entries, 
and $C = U\cdot V$.
Let  $D = \rho^{-1}(C)$, and $E=\rho^{-1}(U)$.
By Lemma  \ref{lemma:key} and Proposition~\ref{prop:commute}, under the mapping $\rho^{-1}$, each successive cell removed from $D$ to obtain $E$ is removed from the last row of the ComT whose length is the column index of the cell being removed, that is the shape of $E$ is $\alpha=row_J(\beta)$, where  $J = \{j_1,\ldots,j_n\}$.
Thus $C = U\cdot V$ is a term in the product $H_\alpha\cdot S_{n}$ from the left hand side.
Moreover, since we are able to uniquely determine $U$ and $V$ from $C$, $C$ appears exactly once on each side of the formula.
This proves the first formula through applying the map $\varepsilon$.

%\bigskip
The proof of the second formula, involving vertical strips, is very much analogous to the first,
making use of the second case of the row bumping lemma.
\end{proof}

\subsection{Transition matrices}\label{sub:transmat}
{}From Theorems~\ref{the:qsasm} and ~\ref{the:qsasf} and the proof of Proposition~\ref{prop:Zbasis} we are able to describe the transition matrices between quasisymmetric Schur functions and monomial or fundamental quasisymmetric functions.

\begin{proposition}\label{prop:q2mmat} Let $A$ be the matrix whose rows  and columns are indexed by $\alpha \vDash n$ ordered by $\btr$ and entry $A_{\alpha\beta}$ is the coefficient of $M_\beta$ in  $\qs _\alpha$. Then $A_{\alpha\beta}$ is the number of ComTs of shape $\alpha$ and weight $\beta$. Furthermore $A_{\alpha\beta} = 0$ if $\alpha \btr \beta$ and $A_{\alpha\alpha}=1$.\end{proposition}

\begin{proof} The first statement follows from Theorem~\ref{the:qsasm}. The second statement follows from the second paragraph of the proof of Proposition~\ref{prop:Zbasis} and the fact that $F_\alpha = \sum _{\alpha \succeq \beta} M_\beta$.
\end{proof}

\begin{proposition}\label{prop:q2fmat} Let $A$ be the matrix whose rows  and columns are indexed by $\alpha \vDash n$ ordered by $\btr$ and entry $A_{\alpha\beta}$ is the coefficient of $F_\beta$ in  $\qs _\alpha$. Then $A_{\alpha\beta}$ is the the number of standard ComTs $T$ of shape $\alpha$ and $\beta(\D(T))=\beta$. Furthermore $A_{\alpha\beta} = 0$ if $\alpha \btr \beta$ and $A_{\alpha\alpha}=1$.\end{proposition}

\begin{proof} The first statement follows from Theorem~\ref{the:qsasf}, while the second statement follows from the second paragraph of the proof of Proposition~\ref{prop:Zbasis}.\end{proof}

The transition matrix from quasisymmetric Schur functions to monomial or fundamental quasisymmetric functions is therefore upper unitriangular by Propositions~\ref{prop:q2mmat} and \ref{prop:q2fmat}.  Consequently, to expand any quasisymmetric function in terms of the quasisymmetric Schur basis, simply invert the appropriate matrix depending on whether the initial quasisymmetric function is given in the monomial or fundamental basis.

Another straightforward application of Proposition~\ref{prop:q2mmat} and Proposition~\ref{prop:q2fmat} yields the following.

\begin{corollary}\label{cor:sisf}
Let $\alpha$ be a composition. Then $\mathcal{S}_\alpha = M_\alpha$ if and only if $\alpha = (1^ f)$. Similarly, $\mathcal{S}_\alpha = F_\alpha$ if and only if $\alpha = (m,1^{e_1},2,1^{e_2},\ldots ,2 ,1^f)$ where $m, f,e_i$ are nonnegative integers such that $m\neq 1$, $f \ge 0$, and $e_i\geq 1$ for all $i$.
\end{corollary}

\begin{example} By Corollary~\ref{cor:sisf} we know that the transition matrix between $\{ \qs _\alpha \} _{\alpha \vDash n}$ and $\{ F _\alpha \} _{\alpha \vDash n}$ is the identity matrix for $n=1,2,3$. For $n=4$ we get
$$\begin{matrix}
{(4)}\\ {(3,1)}\\ {(1,3)}\\ {(2,2)}\\ {(2,1,1)}\\ {(1,2,1)}\\ {(1,1,2)}\\ {(1,1,1,1)}
\end{matrix}\begin{bmatrix}
1&\cdot&\cdot&\cdot&\cdot&\cdot&\cdot&\cdot\\
\cdot&1&\cdot&\cdot&\cdot&\cdot&\cdot&\cdot\\
\cdot&\cdot&1&1&\cdot&\cdot&\cdot&\cdot\\
\cdot&\cdot&\cdot&1&\cdot&1&\cdot&\cdot\\
\cdot&\cdot&\cdot&\cdot&1&\cdot&\cdot&\cdot\\
\cdot&\cdot&\cdot&\cdot&\cdot&1&\cdot&\cdot\\
\cdot&\cdot&\cdot&\cdot&\cdot&\cdot&1&\cdot\\
\cdot&\cdot&\cdot&\cdot&\cdot&\cdot&\cdot&1\\
\end{bmatrix}$$where $\cdot$ denotes 0 and the rows are indexed by quasisymmetric Schur functions. We can hence conclude that our basis differs from those appearing in \cite{BJR, Luoto, Stan}.
\end{example}

\section{Further avenues}\label{sec:conc}

As indicated in the introduction, there are many further avenues to pursue, and in our conclusion we discuss three of them here.

\subsection{A quasisymmetric refinement of the Littlewood-Richardson rule}\label{sub:LR} The Pieri rule generalizes to the celebrated Littlewood-Richardson rule, say \cite[Chapter 7]{ECII}, for expanding the product of two generic Schur functions in terms of Schur functions $$s_\mu s_\nu = \sum _\lambda c^\lambda _{\mu\nu} s_\lambda$$where the Littlewood-Richardson coefficients $c^\lambda _{\mu\nu}$ are positive integers that can be computed combinatorially, given partitions $\lambda, \mu, \nu$. The combinatorial computation requires enumerating all reversetableaux of shape $\lambda /\mu$ and weight $\nu$ subject to one further condition known as the lattice condition. 

Since Theorem~\ref{the:pieriqs} refines the classical Pieri rule, it is natural to ask whether expanding the product of two generic quasisymmetric Schur functions in terms of quasisymmetric Schur functions refines the classical Littlewood-Richardson rule simply. Such a refinement does not presently seem simple, as expanding the product of two generic quasisymmetric Schur functions in terms of quasisymmetric Schur functions often results in negative structure constants. The smallest example is
\begin{align*}
\qs _{(2,1)}\qs _{(2,1)} =\ & \qs _{(4,2)} +  \qs _{(4,1,1)} + 2\qs _{(3,2,1)}+\qs _{(3,1,2)}+ 2\qs _{(2,3,1)}\\
+\ &\qs _{(1,3,2)}+\qs _{(3,1,1,1)}+\qs _{(2,2,2)}+\qs _{(2,2,1,1)}+\qs _{(2,1,2,1)} \\
-\ &\qs _{(1,4,1)}-\qs _{(1,3,1,1)}-\qs _{(1,1,3,1)}-\qs _{(1,2,2,1)}.\end{align*}

However, a product that does naturally refine the classical Littlewood-Richardson rule is the product of a generic \emph{Schur} polynomial with a generic quasisymmetric Schur polynomial expanded in terms of quasisymmetric Schur polynomials. More precisely, in the sequel to this paper \cite{HLMvW2} we prove that
$$s_\lambda (x_1, \ldots , x_n) \qs _\alpha (x_1, \ldots , x_n) = \sum  _\beta C ^\beta _{\alpha\lambda}\qs _\beta (x_1, \ldots , x_n)$$where the $C ^\beta _{\alpha\lambda}$ are positive integers whose computation requires enumerating all ComTs of shape $\beta$ with $\alpha$ removed from the top left corner, with weight the parts of $\lambda$ taken in reverse order, and subject to a lattice-type condition. In addition, we show that similar combinatorial rules exist for the product of a generic Schur polynomial and a Demazure atom   and generic Schur polynomial and a Demazure character   when expanded as a linear combination of Demazure atoms and characters, respectively. Moreover, we recover the classical Littlewood-Richardson rule as a special case of this latter result, when we restrict Demazure characters to Schur polynomials.

\subsection{Skew quasisymmetric functions and duality}\label{sub:skew}
In the classical theory of symmetric functions, say \cite[Chapter 1]{Mac}, there exists the \emph{Hall inner product} $\langle \cdot , \cdot \rangle$, which pairs dual graded bases in the self-dual Hopf algebra $\Lambda$. This inner product reveals that the Schur functions form an orthonormal basis of $\Lambda$, that is
$$\langle s_\lambda , s_\mu \rangle = \delta _{\lambda\mu}$$where $\lambda, \mu$ are partitions and $\delta _{\lambda\mu} = 1$ if $\lambda = \mu$ and $0$ otherwise. Equivalently, the \emph{Cauchy formula} states that
$$\sum _\lambda s_\lambda(x_1, \ldots ) s_\lambda (y_1, \ldots )= \prod _{i,j} \left( 1- x_iy_j\right)^{-1}.$$
One might wonder how such notions extend to the Hopf algebra of quasisymmetric functions $\Q = \bigoplus _{n\geq 0} \Q _n$ and its dual, the algebra of \emph{noncommutative symmetric functions}, $NSym = \bigoplus _{n\geq 0} NSym _n$, introduced in \cite{Gelfand}.

In \cite[Section 6]{Gelfand}, following the work of \cite{Gessel} and \cite{MR}, a pairing between dual graded bases of $\Q$ and $NSym$ was introduced as an analogue to the Hall inner product. This pairing yielded
$$\langle F_\alpha , R_\beta \rangle = \delta _{\alpha\beta}$$where $\alpha, \beta$ are compositions and $R_\beta$ is the noncommutative ribbon Schur function whose commutative image is the ribbon Schur function $r_\beta$. Also, \cite{Gelfand} introduced the equivalent \emph{Cauchy element} in the graded completion of $\bigoplus _{n\geq 0} NSym _n \otimes \Q _n$ as
$$\mathcal{C}:= \sum _\alpha R_\alpha \otimes F_\alpha = \sum _\alpha a_\alpha \otimes b_\alpha$$where $\{ a _\alpha\}$ and $\{ b_\alpha \}$ is any pair of dual graded bases, as a means to describe the dual bases of the various bases of $NSym$.

Conversely, we can ask what can be deduced about the dual basis of quasisymmetric Schur functions $\{ \qs ^* _\alpha \}$ from quasisymmetric Schur functions themselves? For this we need \emph{skew} quasisymmetric Schur functions, and this question is fully addressed in \cite{BLvW}.

\subsection{Quasisymmetric Hall-Littlewood and Macdonald polynomial decompositions}\label{sub:HL}

In view of the fact that the Demazure atoms and characters can be obtained by 
setting $q=t=0$ in various versions of
Macdonald polynomials, a natural question to ask is whether $q$ and/or $t$ parameters
can be inserted in a natural way into the construction of quasisymmetric Schur functions.  
In this section we show 
how the $t$ parameter can  easily be added to some of our constructions, resulting in a decomposition of
the Hall-Littlewood polynomial into quasisymmetric functions.  We contrast this with an 
alternate decomposition obtained from a result in  \cite[Appendix A]{Haglund}, and discuss obstacles 
preventing the insertion of an additional $q$
parameter into our model.  Throughout this section we let $X_n$ denote the ordered sequence of variables
$x_1,\ldots ,x_n$.

Let $\gamma$ be a weak composition into $n$ parts, and $s \in \gamma$, i.e. $s$ a cell or \emph{square} of the diagram of $\gamma$.  Let
$\text{row}(s)$, $\text{col}(s)$, $\text{West}(s)$, and $\text{East}(s)$ 
denote the row containing $s$, the column containing $s$, the square of $\augdg (\gamma)$ in $\text{row}(s)$ 
immediately left of $s$, and the square of $\gamma$ in $\text{row}(s)$ immediately right of $s$ (if it exists), respectively.
Furthermore let $\text{leg}(s)$ be the number of squares in $\text{row}(s)$,  
but to the right of $s$, and $\text{arm}(s)$ the 
number of squares of $\gamma$ in the same column as $s$, below $s$, and in a row not longer than 
$\text{row}(s)$, plus the number of squares of $\augdg (\gamma)$ in the column just left (which may be in the
basement) of $\text{col}(s)$, in a row above
$s$, and also in a row strictly shorter than $\text{row}(s)$.
For a filling $\tau$ of $\gamma$, we let
$\tau (s)$  denote the entry of $\tau$ in $s$.

\begin{example} On the left, the leg lengths, and on the right, the arm lengths, 
for the squares of the augmented diagram (with unmarked basement) $(1,0,3,2,3)$.
$$\tableau{
\ &0\\
\ \\
\ &2&1&0\\
\ &1&0\\
\ &2&1&0}\qquad
\tableau{
\ &0\\
\ \\
\ &4&3&1\\
\ &2&1\\
\ &3&2&1}$$

\label{armleg}
\end{example}

We let $E_{\gamma}^{\prime}(X_n;q,t)$ 
denote the nonsymmetric Macdonald polynomial introduced by
Macdonald in \cite{Mac96} and studied by Cherednik \cite{Che95}, and 
$$
E _{\gamma}(X_n;q,t) = 
E^{\prime} _{\gamma ^{*}}(x_n,\ldots ,x_2,x_1;1/q,1/t)
$$ 
the modified version of the $E^{\prime}$ appearing in work of Marshall \cite{Marshall}, where again $\gamma ^{*} =
(\gamma _{n},\ldots ,\gamma _1)$.
Furthermore let
$\mathcal {E^{\prime}}$ and $\mathcal E$ be the integral forms of the $E^{\prime}$s and 
$E$s, respectively, defined via
\begin{align}
\mathcal {E^{\prime}} _\gamma (X_n;q,t) &= \prod _{s\in \gamma ^{*}}(1-q^{\text{leg}(s)+1}
t^{\text{arm}(s)+1}) \, {E^{\prime}}_{\gamma}(X_n;q,t) \\
\mathcal E _\gamma (X_n;q,t) &= \prod _{s\in \gamma}  (1-q^{\text{leg}(s)+1}
t^{\text{arm}(s)+1}) \, E_{\gamma}(X_n;q,t).
\end{align}
For $\mu$ a partition, we let $P_{\mu}(X_n;q,t)$ denote the symmetric Macdonald
polynomial \cite[Chapter 7]{Mac} and $J_{\mu}(X_n;q,t)$ its integral form
\cite[p. 352]{Mac}, 
\begin{align}
J_{\mu}(X_n;q,t) = \prod_{s\in \mu} (1-q^{\text{leg}(s)}t^{\text{arm}(s)+1})
P_{\mu}(X_n;q,t)
\end{align}
in our notation.  
(They are called integral forms since the coefficients of monomials in them are in $\mathbb Z [q,t]$, while
those in the $E^{\prime}$s, 
$E$s, and $P$s are in $\mathbb Q [q,t]$.)  We note that in $E^{\prime}_{\gamma}$, $E_{\gamma}$ and
$P_{\mu}$, the leading coefficient of $x^{\gamma}$, $x^{\gamma}$, and $x^{\mu}$, respectively, is one where $x^\gamma = x_1^{\gamma _1}x_2^{\gamma _2}\cdots$.
%The $E_{\gamma}$ 
%and $\hat E _{\gamma}$ are normalized so the coefficient of $x^{\gamma}$
%is one.) 

In \cite{HHL} the following combinatorial formula for $\mathcal {E}_{\gamma} (X_n;q,t)$ is 
obtained; 
\begin{align}
\label{master}
\mathcal {E} _{\gamma}(X_n;q,t) = 
\sum_{ \text{non-attacking fillings $\tau$ of $\gamma$} \atop b_i=i}
x^{\tau} q^{\text{maj}(\tau,\gamma)} t^{\text{coinv}(\tau,\gamma)}\\
\notag
\times \prod_{s \in \gamma \atop \tau (s) = \tau (\text {West}(s))} 
(1-q^{\text{leg}(s)+1} t^{\text{arm}(s)+1})
\prod_{s \in \gamma \atop \tau (s) \ne \tau (\text {West}(s))} (1-t),
\end{align}
where $\text{coinv}(\tau,\gamma)$ is the number of triples of the filling 
which are not inversion  triples (i.e. are \emph{coinversion triples}),  
and $\text{maj}(\tau,\gamma)$ is the sum of $\text{leg}(s)+1$, 
over all $s \in \gamma$ where $\tau (\text{West}(s))$ 
is smaller than $\tau (s)$ (i.e. a ``descent").  By $b_i=i$ we mean
the square in the $i$-th row of the basement contains $i$, for $1\le i\le n$. 
As usual, basement squares can be included
in triples.

A nice feature of (\ref{master}) is that if we change the basement to $b_i=n-i+1$,
replace $\gamma$ by $\gamma ^{*}$, 
and sum over non-attacking fillings as above, we get a formula for
$\mathcal {E^{\prime}} _{\gamma}(X_n;q,t)$, while if we sum over non-attacking fillings with
basement $b_i=n+1$ for all $i$, we get
a formula for $J_{\mu}(X_n;q,t)$, where $\mu = \lambda (\gamma)$.
Letting $q=t=0$ in these results give formulas for
Demazure atoms ($\mathcal {E} _{\gamma}(X_n;0,0)$), Demazure characters,
($\mathcal {E^{\prime}} _{\gamma}(X_n;0,0)$) and Schur functions ($J_{\mu}(X_n;0,0)$).

Macdonald obtained an expression for 
$P_{\mu}$ as a linear combination of the $E^{\prime}_{\gamma}$.  
Expressed in terms of the $E$'s, this takes the form \cite{Marshall},
\cite[Eq. (72)]{HHL}
\begin{align}
P_{\mu}(X_n;q,t) = \prod_{s \in \mu} 
(1-q^{\text{leg}(s)+1}t^{\text{arm}(s)}) 
\sum_{\gamma \atop \lambda (\gamma) =\mu}
\frac {E_{\gamma }(X_n;q,t)}{\prod_{s \in \gamma} 
(1-q^{\text{leg}(s)+1}t^{\text{arm}(s)})}.
\end{align}
By setting $q=0$ in this formula
we get
\begin{align}
\label{PtoE}
P_{\mu}(X_n;t) = 
\sum_{\gamma \atop \lambda (\gamma )=\mu}
E_{\gamma }(X_n;0,t)
\end{align}
where $P_{\mu}(X_n;t)=P_{\mu}(X_n;0,t)$ is the Hall-Littlewood polynomial
\cite[p. 208]{Mac}.

It is natural to refer to the function
\begin{align}
E_{\gamma }(X_n;0,t) = 
\mathcal{E} _{\gamma }(X_n;0,t)
\end{align}
as a nonsymmetric Hall-Littlewood polynomial, and we denote this function by $E_{\gamma}(x_1,\ldots,x_n;t)$. 
From (\ref{master}) we have the explicit formula
\begin{align}
\label{masterq}
E _{\gamma}(X_n;t) = 
\sum_{ \text{non-attacking fillings $\tau$ of $\gamma$} \atop b_i=i, \, \,\text{maj}(\tau,\gamma)=0}
x^{\tau} t^{\text{coinv}(\tau,\gamma)}
\prod_{s \in \gamma \atop \tau (s) \ne \tau (\text {West}(s))} (1-t).
\end{align}

For a given composition $\alpha$, let $\mathcal{L}_{\alpha}(X_n;t)$ be the polynomial obtained by 
summing ${E}_{\gamma }(X_n;t)$ over all compositions $\gamma$ for which $\alpha (\gamma) = \alpha$,
\begin{align}
\label{masterql}
\mathcal{L}_{\alpha}(X_n;t) = \sum _{\gamma : \alpha (\gamma) = \alpha}
\sum_{ \text{non-attacking fillings $\tau$ of $\gamma$} \atop b_i=i, \, \,\text{maj}(\tau,\gamma)=0}
x^{\tau} t^{\text{coinv}(\tau,\gamma)}
\prod_{s \in \gamma \atop \tau (s) \ne \tau (\text {West}(s))} (1-t).
\end{align}
Since the quasisymmetric Schur functions $\mathcal{S}_{\alpha}$ are obtained by summing the 
specialization of ${E}_{\gamma}(x;q,t)$ to $q=t=0$ over all compositions which collapse 
to $\alpha$, we have $\mathcal{L}_{\alpha}(X_n;0)=\mathcal{S}_{\alpha}$.
We now show that the $\mathcal {L}_{\alpha}$ are quasisymmetric.
\begin{proposition}
\label{qshl}
The polynomials $\mathcal{L}_{\alpha}$ are quasisymmetric in $x_1,\ldots ,x_n$.
\end{proposition}
\begin{proof}
Note $\mathcal{L}_{\alpha}$ is 
quasisymmetric if and only if the monomial $x_{j_1}^{a_1} x_{j_2}^{a_2} \cdots x_{j_k}^{a_k}$ 
where $j_1 < j_2 < \cdots < j_k $ has the same coefficient as 
$x_{i_1}^{a_1} x_{i_2}^{a_2} \cdots x_{i_k}^{a_k}$ for any other sequence $i_1 < i_2 < \cdots < i_k $.  
We prove this by exhibiting a coinv-preserving bijection between descentless fillings $\sigma$ of a 
weak composition $\gamma$ which collapses to $\alpha$, containing the multiset of 
entries $\{i_1^{a_1},{i_2}^{a_2}, \ldots ,{i_k}^{a_k} \}$, and descentless fillings $\sigma ^{\prime}$ of a 
(possibly different) weak composition which also collapses to $\alpha$, containing the 
multiset of entries $\{j_1^{a_1},{j_2}^{a_2}, \ldots ,{j_k}^{a_k} \}$.  Our bijection will
also preserve the number of squares $s$ of $\gamma$ where $\sigma (s) \ne \sigma (\text{West}(s))$, and
hence will preserve the power of $t$ and $1-t$ multiplying $x^{\sigma}$ in (\ref{masterq}). 

It is straightforward to check that in a descentless, non-attacking filling with basement $b_i=i$, if the entry in 
the first column of a given row is $j$, then the given row must be the $j$-th row.
Let $F$ be such a filling, of a weak composition $\gamma$ with $\alpha (\gamma) = \alpha$, whose 
entries are given by the multiset $\{i_1^{a_1},{i_2}^{a_2}, \ldots ,{i_k}^{a_k} \}$.  
Simply replace each entry $i_s$ by the entry $j_s$ for all $s$ from $1$ to $k$ and slide the 
rows so that the $r^{th}$ row is the row whose first column-entry is $r$.  Note that this preserves 
the order of the nonzero rows since their relative order (given by the entries in the leftmost columns) 
is not affected by the replacement of $i_s$ by $j_s$.  This also implies that the rows remain 
weakly decreasing.  The relative orders of the entries in the triples are preserved, so the 
number of inversion triples is preserved.  Also, squares $s$ where $\sigma (s) \ne \sigma (\text {West}(s))$
are mapped to other such squares, and similarly for squares with 
$\sigma (s) = \sigma (\text {West}(s))$.
To invert this map, simply replace $j_s$ by $i_s$.  
\end{proof}

Proposition \ref{qshl} together with (\ref{PtoE}) imply that 
\begin{align}
\label{Halldecomp}
P_{\mu}(X_n;t) = \sum_{\alpha \atop \lambda (\alpha )=\mu} {\mathcal L}_{\alpha}(X_n;t)
\end{align}
is a decomposition of the Hall-Littlewood polynomial into quasisymmetric functions.
We mention that Hivert \cite{Hivert} has introduced other quasisymmetric 
functions $G_{\alpha}(X_n;t)$ that he calls 
quasisymmetric Hall-Littlewood functions, which he defines via
difference operators.  
Hivert obtains expansions for the $G_{\alpha}$ in terms of the fundamental quasisymmetric
functions $F_{\beta}$, and also in terms of
the monomial quasisymmetric functions $M_{\beta}$, and shows the $G_{\alpha}$ satisfy the 
interesting relations $G_{\alpha}(X_n;0)=F_{\alpha}(X_n)$
and $G_{\alpha}(X_n;1)=M_{\alpha}(X_n)$. 
On the other
hand, ${\mathcal L}_{\alpha}(X_n;0)={\mathcal S}_{\alpha}(X_n)$.  Furthermore, when $t=1$ the only
fillings $\sigma$ defining ${\mathcal L}_{\alpha}$ in (\ref{masterq}) which survive are those for which
there are no squares $s$ with $\sigma (s) \ne \sigma (\text{West}(s))$, i.e. those $\sigma$ which are
constant across rows.  Such $\sigma$ have no coinversions, and it follows that 
the coefficient of $x^{\alpha}$ in ${\mathcal L}_{\alpha}(X_n;1)$ equals $1$, and 
thus ${\mathcal L}_{\alpha}(X_n;1)=M_{\alpha}(X_n)$.
For means of comparison, 
$$G_{13}(X_n;t) = M_{13} + (1-q^2) M_{121} + (1-q^2) M_{112} + (1-2q^2+q^4) M_{1111}$$while

$$\mathcal{L}_{13}(X_n;t) = M_{13}+(1-q)M_{22} + (1-q)M_{211} + (1-q)M_{121} + (2-2q)M_{112} + (2+q)(1-q)^2M_{1111}$$where we drop the brackets around the compositions for brevity.

In \cite{Haglund} an explicit decomposition of $J_{\mu}(X_n;q,t)$ 
into the $F_{\alpha}$ is obtained.  Since this formula has not appeared in a journal article
before, we include a detailed description of it here, and contrast the $q=0$ case of it
with  
the decomposition of $P_{\mu}(X_n;t)$ into the 
quasisymmetric functions $\mathcal L_{\alpha}$ above.
Consider a standard filling  $\tau$ of $\mu$, with basement
$(n+1,\ldots ,n+1)$.  Such a filling is automatically non-attacking, and $\tau$ can be identified
with the permutation obtained by reading in the entries of $\tau$, from top to bottom within columns,
starting with the rightmost column and working right to left. 
Given a triple of $\tau$ (neccessarily of type A since $\mu$ is a partition) which 
doesn't involve any basement squares, we call the square containing the middle
of the three entries (i.e. neither the largest nor the smallest) 
the ``base" of the triple.  If
the triple involves a basement square, we call the square containing the smallest of the
three entries the base of the triple. 

\begin{example} 
A standard filling of $(3,3,1)$, with basement $(8,8,8)$. 
$$\tableau{
8&{\bf 5}&{\bf 6}&{\bf 1}\\
8&{\bf 2}&{\bf 7}&{\bf 4}\\
8&{\bf 3}}$$ For   the filling above,
the base square of the triple consisting of entries $5,6,7$ contains the $6$, the triple with
entries $1,4,6$ has base containing the $4$, the triple with entries $2,3,8$ has base containing
the $2$, and the base
square of the triple consisting of entries $3,5,8$ contains the $3$.
\label{base}
\end{example}

Let $\text{coinv}_s(\tau,\mu)$ be the number of coinversion triples, and
$\text{inv}_s(\tau,\mu)$ the number of inversion triples, where $s$ is the base
square.  
Also, if $\tau (\text{East}(s)) > \tau (s)$ (so there is a descent at $\text{East}(s)$), let
$\text{maj}_s (\tau,\mu) = \text{leg}(s)$, else set 
$\text{maj}_s (\tau,\mu) = 0$.  And, 
if $\tau (\text{West}(s)) \ge \tau (s)$ (so there is no descent at $s$), let
$\text{nondes}_s (\tau,\mu) = \text{leg}(s)+1$, else set 
$\text{nondes}_s (\tau,\mu) = 0$.  Note that 
$\text{coinv}(\tau,\mu)=\sum_{s \in \mu} \text{coinv}_s(\tau,\mu)$, with similar statements for
$\text{inv}(\tau, \mu)$ and $\text{maj}(\tau, \mu)$.   Then we have \cite[p. 133]{Haglund}
\begin{align}
\label{JmuQ}
J_{\mu}(X_n;q,t) =
\sum_{\tau \in S_n }
F_{\beta (\{i : \tau ^{-1}(i)>\tau ^{-1}(i+1)\})}(X_n)
\prod _{ s \in \mu}
(q^{\text{inv}_s(\tau ,\mu)}t^{\text{nondes}_s(\tau ,\mu)}
-q^{\text{coinv}_s(\tau ,\mu)}t^{1+\text{maj}_s(\tau ,\mu)}),
\end{align}
which gives an expansion of $J_{\mu}$ into Gessel's fundamental quasi-symmetric
functions.  
It is also shown in \cite{Haglund} that if $\tau$ is such that
some entry $j$ occurs in a column to the right of the $j$-th column, then the
factor
\begin{align}
\label{Qterm}
\prod_{s \in \mu}
(q^{\text{inv}_s(\tau,\mu)}t^{\text{nondes}_s(\tau,\mu)}
-q^{\text{coinv}_s(\tau,\mu)}t^{1+\text{maj}_s(\tau,\mu)})
\end{align}
in (\ref{JmuQ}) is zero.

By letting $q=0$ in (\ref{JmuQ}), we get a decomposition of the integral form Hall-Littlewood
polynomial ($Q_{\mu}(X_n;t)$ in the notation of \cite[p. 210]{Mac}) into fundamental
quasisymmetric functions.  It is complicated, though, to work with the set of permutations over
which (\ref{Qterm}) does not vanish when $q=0$, i.e. the set where every
square $s$ of $\mu$
satisfies either $\text{coinv}_s(\tau,\mu)=0$ or $\text{inv}_s(\tau,\mu)=0$, or both.  Also,
the formula for ${\mathcal L}_{\alpha}$ as a sum over non-attacking fillings 
is a positive formula in the sense that
each coefficient of a monomial is a sum of terms of the form $t^* (1-t)^*$ for
nonnegative integers $*$, while the monomial coefficients in the
$q=0$ case of formula (\ref{JmuQ}) could involve terms of the form $\pm t^*\prod (1-t^*)$.  We also mention
that we need to divide the $q=0$ case of (\ref{JmuQ}) by the product
\begin{align}
\prod_{s\in \mu \atop \text{leg}(s)=0} (1-t^{\text{arm}(s)+1})
\end{align}
to convert from the integral form $J_{\mu}(X_n;0,t)$ to $P_{\mu}(X_n;t)$, and once we do the coefficient
of a given monomial in the $x$'s is a rational function in $t$, not clearly a polynomial.

One would naturally hope to insert a $q$ parameter into the construction of
$\mathcal {L}_{\alpha }(X_n;t)$, and end up with a decomposition of
$J_{\mu}(X_n;q,t)$ into quasisymmetric functions, where the quasisymmetric
extension of the $\mathcal {L}_{\alpha}(X_n;t)$ is a positive sum in the sense of the above
paragraph, along the lines of the formula (\ref{master}).  The problem is that the bijective map from the
proof of Proposition \ref{qshl} does not apply as is to fillings with descents just right of basement
squares.  Thus at this time the authors do not see how to extend the construction of the $\mathcal{L}_{\alpha}$ in an
elegant way to include the $q$ parameter.   Another interesting question we leave for future research is how to
decompose the $\mathcal{L}_{\alpha}$ into fundamental quasisymmetric functions.

\section{Acknowledgements} The authors would like to thank Ole Warnaar and the referees for helpful comments and suggestions.


\begin{thebibliography}{130}
\bibliographystyle{abbrv}

\bibitem{Aguiar} M. Aguiar, N. Bergeron and F. Sottile, Combinatorial Hopf algebras and generalized Dehn-Sommerville equations, \emph{Compos. Math.} \textbf{142} (2006),  1-30.

\bibitem{BergBousDul} F. Bergeron, M. Bousquet-Melou and S. Dulucq, Standard paths in the composition poset,  \emph{Ann. Sci. Math. Qu\'{e}bec}  \textbf{19} (1995), 139--151.

\bibitem{BMSvW} N. Bergeron, S. Mykytiuk, F. Sottile and S. van Willigenburg, Non-commutative Pieri operators on posets, \emph{J. Combin. Theory Ser. A} \textbf{91}  (2000), 84--110.

\bibitem{BH} B. Bernevig and F.  Haldane, Model fractional quantum Hall states and Jack polynomials,
{\it Phys. Rev. Lett.} \textbf{100}  (2008), 246802.

\bibitem{BLvW} C. Bessenrodt, K. Luoto and S. van Willigenburg, Skew quasisymmetric Schur functions, in preparation.

\bibitem{BHvW} L. Billera, S. Hsiao and S. van Willigenburg, Peak quasisymmetric functions and Eulerian enumeration, \emph{Adv. Math.} \textbf{176}  (2003), 248--276.

\bibitem{BJR} L. Billera, N. Jia and V. Reiner, A quasisymmetric function for matroids,
\emph{European J. Combin.}, to appear.

\bibitem{HDL}
L. Billera, H. Thomas and S. van Willigenburg,
Decomposable compositions, symmetric quasisymmetric functions and 
equality of ribbon Schur functions, {\it Adv.  Math.} {\bf 204} (2006), 204--240.

\bibitem{Cherednik} I. Cherednik, 
Double affine Hecke algebras and Macdonald's conjectures,
\emph{Ann. Math.} \textbf{141} (1995), 191--216.

\bibitem{Che95} I. Cherednik, 
Nonsymmetric Macdonald polynomials, \emph{Int. Math. Res. Not.}
\textbf{10} (1995), 483--515.

\bibitem{Gelfand} I. Gel'fand, D. Krob, A. Lascoux, B. Leclerc,
V. Retakh and J.-Y. Thibon, Noncommutative symmetric functions,
{\it Adv.  Math.} {\bf 112} (1995), 218--348.

\bibitem{Gessel}  I. Gessel, 
Multipartite $P$-partitions and inner products of skew Schur functions,
combinatorics and algebra
(Boulder, Colo., 1983)  289--317, Contemp. Math., 34, Amer. Math. Soc.,
Providence, RI, 1984.

\bibitem{Haglund} J. Haglund,
     \emph{The $q$,$t$-Catalan numbers and the space of diagonal
              harmonics
      (With an appendix on the combinatorics of Macdonald
              polynomials)},
    American Mathematical Society, University Lecture Series, Volume 41,
Providence, RI, USA, 2008.

\bibitem{HHL2} J. Haglund, M. Haiman and N. Loehr, A combinatorial formula for Macdonald polynomials, 
\emph{J. Amer. Math. Soc.}  {\bf 18}  (2005), 735--761.

\bibitem{HHL}
J. Haglund, M. Haiman and N. Loehr,
  A combinatorial formula for nonsymmetric {M}acdonald polynomials,
  {\it Amer. J. of Math.}  {\bf 103}  (2008),  359--383.
  
  \bibitem{HLMvW2} J. Haglund, K. Luoto, S. Mason and S. van Willigenburg, Refinements of the Littlewood-Richardson rule, ArXiv: 0908.3540.
  
  \bibitem{Hivert} F. Hivert, Hecke algebras, difference operators and quasi-symmetric functions, \emph{Adv. Math.} {\bf 155} (2000), 181--238. 
  
\bibitem{Ion} B. Ion, Nonsymmetric Macdonald polynomials and Demazure characters,
  \emph{Duke Math. J.} {\bf 116} (2003), 299--318.   
  
\bibitem{Ion2} B. Ion, Standard bases for affine parabolic modules and nonsymmetric Macdonald polynomials,
 \emph{ J. Algebra} {\bf 319} (2008),  3480--3517.   

\bibitem{Kad88}
K.  Kadell,
A proof of some analogues of Selberg's integral for $k=1$,
{\it SIAM J. Math. Anal.} {\bf 19} (1988), 944--968.

\bibitem{KrobThibon} D. Krob and J.-Y. Thibon, Noncommutative symmetric functions V: a degenerate version of $U_q(gl_N)$, \emph{Int. J. Alg. Comput.} {\bf 9} (1999), 405--430.

\bibitem{LS} A. Lascoux and M.-P. Sch\"{u}tzenberger, Keys and standard bases, 125--144, Tableaux and Invariant Theory, IMA Volumes in Mathematics and its Applications, vol. 19, 1990.

\bibitem{Luoto} K. Luoto, A matroid-friendly basis for the quasisymmetric functions,
{\it J. Combin. Theory Ser. A}
\textbf{115} (2008), {777--798}.

\bibitem{Mac88} I. Macdonald, A new class of symmetric functions, \emph{S\'{e}m. Lothar. Combin.}
  \textbf{372} (1988), 131--171.

\bibitem{Mac} 
I. Macdonald, {\em Symmetric functions and {H}all polynomials, 2nd Edition}, 
 Oxford University
Press, New York, USA, 1995.

\bibitem{Mac96} I. Macdonald, Affine {H}ecke algebras and orthogonal polynomials, \emph{Ast\'erisque} {\bf 237} (1996),
189--207, S\'eminaire Bourbaki $1994/95$, Exp. no. $797$.

\bibitem{MacM} P. MacMahon, Combinatory Analysis, Dover Publications, 2004.

\bibitem{MR} C. Malvenuto and C. Reutenauer, Duality between quasi-symmetric
functions and the Solomon descent algebra, {\it J.  Algebra} {\bf 177}
(1995), 967--982.

\bibitem{Marshall} D. Marshall, Symmetric and nonsymmetric Macdonald polynomials,
\emph{Ann. Comb.} {\bf 3} (1999), 385--415.

\bibitem{Mason1}
S. Mason, A decomposition of {S}chur functions and an analogue of the {R}obinson-{S}chensted-{K}nuth algorithm, \emph{S\'{e}m. Lothar. Combin.} \textbf{57} (2008), B57e.

\bibitem{Mason2} S. Mason, An explicit construction of Type $A$ Demazure atoms,	\emph{J. Algebraic Combin.} {\bf 29} (2009), 295--313.

\bibitem{ECII}
R. Stanley,
\emph{Enumerative Combinatorics, Volume~2},
Cambridge University Press, Cambridge, UK, 1999. 

\bibitem{Stan} R. Stanley, The descent set and connectivity set of a permutation,
\emph{J. Integer Seq.} \textbf{8} (2005),  05.3.8. 

\bibitem{SvW} S. van Willigenburg, Equality of Schur and skew Schur functions, \emph{Ann. Comb.} {\bf 9} (2005), 355--362.



\end{thebibliography}
\end{document}